\documentclass[preprint,11pt,authoryear]{elsarticle}

\usepackage{amssymb}
\usepackage{amsmath}

\usepackage{setspace}

\usepackage{amsfonts}
\usepackage{algorithmic}
\usepackage{algorithm}
\usepackage{array}
\usepackage[caption=false,font=footnotesize]{subfig}
\usepackage{tabularx}
\usepackage{textcomp}
\usepackage{stfloats}
\usepackage{url}
\usepackage{verbatim}
\usepackage{graphicx}
\usepackage{cite}
\usepackage{accents}
\usepackage{booktabs}
\usepackage{enumitem}
\usepackage{multicol}
\usepackage{multirow}
\usepackage[dvipsnames]{xcolor}
\newlist{abbrv}{itemize}{1}
\setlist[abbrv,1]{label=,labelwidth=1.0in,align=parleft,itemsep=0.1\baselineskip,leftmargin=!}

\usepackage{pifont}
\newcommand{\cmark}{\ding{51}}%
\newcommand{\xmark}{\ding{55}}%

\usepackage{geometry}
\usepackage{tikz}
\usetikzlibrary{shapes.geometric}

\usepackage{balance}
\usepackage{todonotes}

\usepackage{lipsum}
\makeatletter
\def\ps@pprintTitle{%
 \let\@oddhead\@empty
 \let\@evenhead\@empty
 \def\@oddfoot{}%
 \let\@evenfoot\@oddfoot}
\makeatother

\begin{document}

\begin{frontmatter}



\title{
 On the Viability of Stochastic Economic Dispatch \\ for Real-Time Energy Market Clearing
}


\author{Haoruo Zhao} 
\author{Mathieu Tanneau}
\author{Pascal Van Hentenryck}

\affiliation{
            organization={H. Milton Stewart School of Industrial and Systems Engineering, Georgia Institute of Technology}, 
            addressline={755 Ferst Dr NW}, 
            city={Atlanta},
            postcode={30318}, 
            state={GA},
            country={USA}
            }

\begin{abstract}

Over the past decade, the rapid adoption of intermittent renewable energy sources (RES), especially wind and solar generation, has posed challenges in managing real-time uncertainty and variability. In the U.S., Independent System Operators (ISOs) solve a security-constrained economic dispatch (SCED) every five minutes to clear real-time electricity markets, co-optimizing energy dispatch and reserve to minimize costs while meeting physical and reliability constraints. All SCED formulations in the U.S. are deterministic and mostly consider a single time period, limiting their effectiveness in managing real-time operational uncertainty from RES intermittency. This limitation is highlighted by the recent introduction of multiple short-term ramping products in U.S. markets, aiming to bridge the gap between deterministic and stochastic SCED formulations.
While stochastic formulations address uncertainty in a unified and endogenous manner, their adoption has been hindered by high computational costs and, to a lesser extent, the availability of probabilistic forecasts.
This paper revisits these concerns and demonstrates that stochastic economic dispatch is now a viable technology for real-time market clearing. It introduces the stochastic look-ahead dispatch (SLAD) formulation for real-time market clearing and presents an accelerated Benders' decomposition to solve it efficiently. Extensive experiments on a real, industry-sized transmission grid demonstrate the computational scalability of the proposed approach, with SLAD instances being solved in under 5 minutes. Furthermore, results show that SLAD provides more than 50\% additional savings compared to flexiramp products and is more robust to the forecasting methodology. Therefore, SLAD is a promising approach for uncertainty management in real-time electricity markets.

\end{abstract}


\end{frontmatter}



\newcommand{\ubar}[1]{\underaccent{\bar}{#1}}
\newcommand{\RMP}[1]{\text{RMP}({#1})}


\newcommand{\pg}{\mathbf{pg}}
\newcommand{\res}{\mathbf{res}}
\newcommand{\reg}{\mathbf{r}^{\text{reg}}}
\newcommand{\spin}{\mathbf{r}^{\text{spin}}}
\newcommand{\suppOn}{\mathbf{r}^{\text{s,on}}}
\newcommand{\suppOff}{\mathbf{r}^{\text{s,off}}}
\newcommand{\rcup}{\mathbf{rc}^{\uparrow}}
\newcommand{\rcdn}{\mathbf{rc}^{\downarrow}}
\newcommand{\rcstr}{\mathbf{rc}^{\text{str}}}

\newcommand{\pf}{\mathbf{pf}}

\newcommand{\dpp}{\mathbf{\boldsymbol{\delta} p}^{+}}
\newcommand{\dpm}{\mathbf{\boldsymbol{\delta} p}^{-}}
\newcommand{\dr}{\mathbf{\boldsymbol{\delta} r}}
\newcommand{\dres}{\mathbf{\boldsymbol{\delta} res}}
\newcommand{\dreg}{\dr^{\text{reg}}}
\newcommand{\dregspin}{\dr^{\text{rspin}}}
\newcommand{\dopres}{\dr^{\text{op}}}
\newcommand{\drcup}{\mathbf{\boldsymbol{\delta} rc}^{\uparrow}}
\newcommand{\drcdn}{\mathbf{\boldsymbol{\delta} rc}^{\downarrow}}
\newcommand{\drcstr}{\mathbf{\boldsymbol{\delta} rc}^{\text{str}}}
\newcommand{\df}{\mathbf{\boldsymbol{\delta} f}}
\newcommand{\p}{\mathbf{p}}

\newcommand{\x}{\mathbf{x}}
\newcommand{\y}{\mathbf{y}}

\newcommand{\rrup}{\text{r}^{\uparrow}}  
\newcommand{\rrdn}{\text{r}^{\downarrow}}  
\newcommand{\CF}{\text{CF}}  
\newcommand{\RF}{\text{RF}}  
\newcommand{\ITF}{\text{IF}}  
\newcommand{\regAvail}{\text{RA}^{\text{reg}}}
\newcommand{\spinAvail}{\text{RA}^{\text{spin}}}
\newcommand{\suppOnAvail}{\text{RA}^{\text{s,on}}}
\newcommand{\suppOffAvail}{\text{RA}^{\text{s,off}}}
\newcommand{\pfmin}{\ubar{\text{pf}}}
\newcommand{\pfmax}{\bar{\text{pf}}}
\newcommand{\pgmin}{\ubar{\text{pg}}}
\newcommand{\pgmax}{\bar{\text{pg}}}
\newcommand{\pgstep}{\bar{\text{pg}}}
\newcommand{\dt}{\Delta_{t}}
\newcommand{\rt}{\Gamma_{t}}
\newcommand{\pd}{\text{pd}}  
\newcommand{\PTDF}{\Phi}
\newcommand{\regmax}{\bar{\text{r}}^{\text{reg}}}
\newcommand{\spinmax}{\bar{\text{r}}^{\text{spin}}}
\newcommand{\suppOnmax}{\bar{\text{r}}^{\text{s,on}}}
\newcommand{\suppOffmax}{\bar{\text{r}}^{\text{s,off}}}
\newcommand{\rcupmax}{\bar{\text{rc}}^{\uparrow}}
\newcommand{\rcdnmax}{\bar{\text{rc}}^{\downarrow}}
\newcommand{\maxOfflineResponse}{\text{MOR}}

\newcommand{\resReq}{\text{rr}^{\text{res}}}
\newcommand{\regReq}{\text{rr}^{\text{reg}}}
\newcommand{\regSpinReq}{\text{rr}^{\text{rspin}}}
\newcommand{\opResReq}{\text{rr}^{\text{op}}}
\newcommand{\rcupReq}{\text{rr}^{\rcup}}
\newcommand{\rcdnReq}{\text{rr}^{\rcdn}}
\newcommand{\rcstrReq}{\text{rr}^{\rcstr}}

\newcommand{\noload}{\text{c}^{\text{NL}}}
\newcommand{\costPg}{\text{c}^{\text{pg}}}
\newcommand{\costRes}{\text{c}^{\text{res}}}
\newcommand{\costRegRes}{\text{c}^{\text{reg}}}
\newcommand{\costSpinRes}{\text{c}^{\text{spin}}}
\newcommand{\costSuppOnRes}{\text{c}^{\text{s,on}}}
\newcommand{\costSuppOffRes}{\text{c}^{\text{s,off}}}

\newcommand{\penaltyPowerShortage}{\pi^{-}}
\newcommand{\penaltyPowerSurplus}{\pi^{+}}
\newcommand{\penaltyRegRes}{\pi^{\text{reg}}}
\newcommand{\penaltyRegSpinRes}{\pi^{\text{rspin}}}
\newcommand{\penaltyOpRes}{\pi^{\text{op}}}
\newcommand{\penaltyRCUp}{\pi^{\text{rc}^{\uparrow}}}
\newcommand{\penaltyRCDn}{\pi^{\text{rc}^{\downarrow}}}
\newcommand{\penaltyRCStr}{\pi^{\text{rc}^{\text{str}}}}
\newcommand{\penaltyFlow}{{\pi^{\text{f}}}}

\section{Introduction}
\label{sec:intro}

The Security-Constrained Economic Dispatch (SCED) is a fundamental optimization problem in power systems operations.
In the US, Independent System Operators (ISOs) clear the real-time electricity market by solving a SCED every few minutes.
The SCED outputs energy and reserve dispatch for each generator, as well as energy and reserve prices \citep{Ma2009}.

Over the past few years, the growing penetration of renewable generation and distributed energy resources (DERs) has caused an increase in operational uncertainty and short-term flexibility requirements.
Inadequate flexibility leads to price spikes, expensive emergency actions, and, in the worst case, load shedding events.
As a result, ISOs introduced flexible ramping products --also known as \emph{flexiramp}-- in electricity markets, with the goal of mitigating short-term uncertainty.
For instance, the Midcontinent ISO (MISO) implemented a 10-minute up/down ramping product in 2016 \citep{wang_hansen_chatterjee_merring_li_li_2016}, and introduced a 30-minute ramping product in December 2021 \citep{Chen2023_MISORPDesign, MISO2021_ShortTermReserve}.
Similarly, the California ISO (CAISO) market includes a 5-minute ramping product \citep{CAISO_BPM}.

While it provides reliability benefits, {\em flexiramp} is outperformed by stochastic optimization models, as convincingly demonstrated in \citep{WangHobbs2013_Flexiramp}.
Stochastic look-ahead dispatch (SLAD) considers the uncertainty of future events and seeks to find the optimal dispatch decisions that maximize the system performance, such as minimizing costs and maximizing reliability.
Nevertheless, no US ISO currently uses stochastic optimization in production, in a large part because of the increased computational cost \citep{WangHobbs2013_Flexiramp}.
Furthermore, to the best of the authors' knowledge, {\em there has been no study that quantifies the potential benefits of stochastic economic dispatch formulations on real, industrial systems compared to existing deterministic approaches with ramping products.}

This paper addresses these gaps by proposing SLAD as a stochastic extension of existing deterministic formulations for real-time electricity markets.
The benefits of SLAD are demonstrated on a realistic, large-scale power system, through simulations that replicate the operations of an ISO.
Moreover, the paper presents a comprehensive comparison of several real-time dispatch approaches, ranging from simple myopic economic dispatch, to multi-period and stochastic dispatches.

\subsection*{Contributions and outline}
\label{sec:intro:contribution}

Building upon the challenges identified by Hobbs, this study addresses the management of increasingly volatile net loads driven by the growth in renewable energy. 
Utilizing stochastic look-ahead dispatch SLAD, the paper explores the benefits of multi-period formulations and ramping products within a practical, industrial context.
The paper's key contributions are the following:

\begin{enumerate}
    \item The formulation of a SLAD model and an efficient solution methodology to solve it.
        The formulation integrates energy and reserve variables.
        The solution algorithm uses parallel computing and accelerated Benders decomposition for enhanced computational efficiency.
    \item Simulations on a large-scale realistic network confirm SLAD's superior performance. Comprehensive comparisons demonstrate its advantages over traditional deterministic real-time dispatch models with ramping products, under realistic forecasts and scenarios.
    \item A detailed case study and sensitivity analysis reveal why SLAD outperforms existing approaches and to show how different scenarios and look-ahead horizons affect its effectiveness, providing insights for system planning and operations.
\end{enumerate}

\noindent
The rest of the paper is organized as follows.
Section \ref{sec:literature} surveys the relevant literature.
Section \ref{sec:formulation} introduces the economic dispatch formulations considered in the paper.
Section \ref{sec:method} presents the proposed accelerated Benders decomposition algorithm for solving SLAD.
Section \ref{sec:example} illustrates the different formulations on a small example.
Section \ref{sec:experiment_RTE} presents the experiment setting on a real-life system, Section \ref{sec:results} presents numerical results, and Section \ref{sec:sensitivity} discusses the
sensitivity of each formulation.
Section \ref{sec:conclusion} identifies future research directions and concludes the paper.

\section{Related Work}
\label{sec:literature}

\subsection{Deterministic Formulation}

\paragraph{Security-Constrained Economic Dispatch}

The security-constrained economic dispatch (SCED) is the mathematical formulation that underlies virtually all real-time electricity markets in the US \citep{CAISO_BPM,NYISO_BPM,MISO_BPM002_D,ERCOT_BPM}.
The SCED formulations used in industry are linear programming (LP) problems that decide each generator's energy and reserve dispatch so as to minimize total operating costs, subject to physical, engineering and market constraints \citep{Ma2009}.
In its simplest form, SCED is a single-period, deterministic LP problem with tens of thousands of variables and constraints.
It is typically solved in a few seconds using state-of-the-art LP solvers.

\paragraph{SCED with Flexiramps}

The main limitation of a single-period SCED formulation is that it does not consider the uncertainty and volatility of the net load.
Therefore, to manage operational uncertainty in real-time markets, US ISOs are increasingly using short-term ramping products, also known as Flexible Ramping Products (FRPs) or \emph{flexiramp}.
For instance, MISO's real-time market includes a 10-minute and a 30-minute ramping product \citep{Chen2023_MISORPDesign,MISO_BPM002_D}.
Similarly, CAISO considers a 5-minute ramping product in its real-time market \citep{CAISO_BPM}. Also, similar ideas are proposed for unit commitment problems \citep{Khodayar2016_MultiplePeriodRamping}.
Nevertheless, previous studies on the benefits of these FRPs have observed that single-period formulations may still fail to account for future ramping needs, ultimately resulting in sub-optimal decisions \citep{WangHobbs2013_Flexiramp}.
While FRPs can improve system flexibility, their design and implementation present challenges such as inaccurate forecasting, uncertainty estimation, and lack of structured pricing mechanisms \citep{SREEKUMAR2022112429}.
Therefore, the authors recommend a careful ramp product market design and advocate a shift towards stochastic dispatch and commitment models.
However, due to differences and incompatibilities in market design, an alternative method is proposed for European markets to approximate stochastic dispatch by determining reserve requirements for zonal operating reserves \citep{Dvorkin2019_SettingReserveRequirements}.

\paragraph{Look-Ahead Dispatch}

To extend SCED from a single-period to a multi-period formulation, look-ahead dispatch (LAD) is proposed in \citep{Xie2009_LAD_MPC}.
The main motivation for LAD is that the longer horizon allows to take into account future operating conditions, thereby achieving better economic performance in the long run.
While the use of a look-ahead window is well established in ISOs' intra-day reliability unit commitment process \citep{Chen2012_MISO_LAC,Hui2009_ReliabilityUnitCommitmentERCOT,Ott2010_ComputingRequirementsPJMMarket}, the use of LAD formulations in real-time markets has received less attention.

A LAD formulation for the ERCOT nodal market was been proposed in \citep{Xie2011_LookAheadDispatchERCOT}. 
Building on this work, a study on the ERCOT system was conducted wherein a LAD with a one-hour look-ahead horizon is considered \citep{Hui2013_ERCOT_LAD}.
The paper raises the question of ``\emph{how far should LAD look ahead?}", and highlights that forecast errors and unforeseen events can create misleading price signals, thereby resulting in sub-optimal decisions.
The observation that LAD is sensitive to forecast errors is echoed in \citep{Choi2017_DataPerturbationLAD}, which evaluates the sensitivity of LAD to changes in its input data.
In particular, the authors report that LAD may exhibit worse performance than a single-period SCED when errors in the input data are too large.
More recently, a reduced formulation for a LAD model that considers the N-k contingency criterion has been proposed, with results reported on systems with up to 1354 buses \citep{VARAWALA2023}.
However, the reported computing times exceed 10 minutes for the larger case, which is too high for real-time market clearing (which occurs every 5 minutes).

In the US, both CAISO and NYISO use a LAD formulation in their real-time markets \citep{CAISO_BPM,NYISO_BPM}.
Namely, CAISO considers a one-hour look-ahead window with 5-minute time intervals, and NYISO considers a one-hour look-ahead window with a mixed 5-minute and 15-minute time intervals.
Dispatch decisions in the first interval of LAD are binding, and dispatch decisions from subsequent intervals are advisory.

\subsection{Stochastic Formulation}

The reader is referred to \citep{Roald2023_MethodsReview} for a detailed overview of power systems optimization under uncertainty, which includes stochastic, robust, distributionally robust and chance-constrained optimization.
The latter two approaches are computationally heavy, and have not been demonstrated to scale to industry-size instances.
While several studies have used robust optimization (RO) for economic dispatch \citep{Li2015_AdujstableROEconomicDispatch,Lorca2015_AdaptiveROEconomicDispatch,YILDIRAN2023890,Zhao2015_MISORobustOptimization}, most consider small systems with only a few hundred buses, which is two orders of magnitude smaller than real-life instances.
Furthermore, the tractability of RO heavily depends on the choice of uncertainty set: the resulting RO problem may be, e.g., a linear, second-order cone, mixed-integer, or semi-definite positive optimization problem.
The design of data-driven, computationally tractable uncertainty sets in high-dimensional settings is an active area of research.

In contrast, stochastic optimization takes as input a probabilistic forecast in the form of scenarios, which can be generated using a variety of models \citep{mashlakov2021assessing,Messner2020_EvaluationWindPowerForecast,werho2021scenario}.
Importantly, the structure of stochastic optimization models, and the corresponding solution algorithms, is agnostic to the nature of the probabilistic forecasting model.
Therefore, the paper focuses on stochastic optimization approaches, and assumes that a probabilistic forecasting model is available.
The training of probabilistic forecast models is beyond the paper's scope; the reader is referred to \citep{mashlakov2021assessing} for a recent overview of available tools.

Stochastic look-ahead dispatch (SLAD) (see Section \ref{sec:formulation:SLAD}) is a stochastic version of LAD to hedge against forecast uncertainty.
Previous work has shown that, compared to deterministic ED formulations, SLAD can improve the reliability and efficiency.
In \citep{Reynolds2020Scenario}, SLAD is evaluated on a 51-bus system, demonstrating more efficient dispatch decisions compared to SCED. 
A multistage SLAD model with renewable production uncertainty and storage is proposed in \citep{Papavasiliou2018SDDP}, and solved using Stochastic Dual Dynamic Programming (SDDP).
The approach is evaluated on a 228-bus system representing the German grid, with results showing around 1\% savings compared to SCED . 
Similarly, a study on SLAD with a flexible ramping product demonstrates its effectiveness and resilience in handling worst-case scenarios on a 5-bus system \citep{Zhang2015_SLADWithFlexiramp}. 
Further comparisons of SLAD with SCED and LAD on a 5889-bus system in \citep{gu2016stochastic} indicate that SLAD provides greater cost savings, especially during intervals with high economic risks. 
However, this formulation does not account for reserve variables, which are crucial in U.S. real-time markets \citep{gu2016stochastic}.

One of the primary challenges in stochastic optimization is accurately representing the uncertainty.
To address this challenge, the sample average approximation method is often used.
However, this method typically requires a large number of samples to estimate the expected cost accurately.
As a result, the problem size increases substantially, making it challenging to solve the optimization problem efficiently.
Benders decomposition is a well-suited method for solving stochastic optimization problems; it leverages the inherent structure of these formulations \citep{conejo2006decomposition}.
SDDP is a method used for solving multi-stage stochastic problems \citep{PEREIRA1989161,Pereira1991359} by combining cutting-plane methods with a Monte Carlo simulation.
Progressive hedging \citep{boyd2011distributed} is another technique used to address two-stage stochastic problems.
This method involves decomposing the problem by considering an augmented Lagrangian of the original problem.

Finally, the paper focuses on addressing the scalability challenges of SLAD for industry-size problems, and evaluates the potential operational benefits of SLAD on a real industrial system.
The design of SLAD-compatible pricing schemes, while important for the adoption of SLAD in real-time markets, is beyond the scope of the paper.
The reader is referred to \citep{chen2020pricing,guo2021pricing,Hua2019Price,Tan2022LMPReview,ZhaoEtAl2020_MultiPeriodMarket} for recent studies on such pricing schemes.

\section{Problem Formulation}
\label{sec:formulation}

This section presents the single-period (SCED), single-period with ramping products (SCED+RP), 
deterministic multi-period (LAD), and stochastic multi-period (SLAD) economic dispatch formulations considered in the paper.
These formulations are based on MISO's real-time market formulation \citep{MISO_BPM002_D}.
Table \ref{tab:nomenclature} presents a nomenclature of the notations used in the paper.
Table \ref{tab:abstract_formulations} outlines the number of time steps, the number of scenarios, 
and whether flexible ramping products (FRPs) are included in each problem formulation.

\begin{table}[!t]
    \centering
    \caption{Nomenclature}
    \label{tab:nomenclature}
    \small
    \begin{tabularx}{\textwidth}{lX|lX}
        \toprule
        \multicolumn{4}{l}{\textbf{Sets and Indices}} \\
        \midrule
        $i \in \mathcal{N}$ & Buses & $e \in \mathcal{E}$ & Branches \\
        $g \in \mathcal{G}$ & Generators & $t \in \mathcal{T}$ & Time steps \\
        $\mathcal{G}_{i}$ & Generators at bus $i$ & $s \in \mathcal{S}$ & Scenarios \\
        \midrule
        \multicolumn{4}{l}{\textbf{Variables}} \\
        \midrule
        $\p_{i,t,s}$ & Net power injection at bus $i$ at time $t$ in scenario $s$ (MW) & $\rcup_{g, t, s}, \rcdn_{g, t, s}$ & Ramping-up and -down capability of generator $g$ at time $t$ in scenario $s$ (MW) \\
        $\pg_{g, t, s}$ & Energy dispatch of generator $g$ at time $t$ in scenario $s$ (MW) & $\drcup_{t, s}, \drcdn_{t, s}$ & Ramping-up and -down capability shortage at time $t$ in scenario $s$ (MW) \\
        $\res_{g,t,s}$ & Reserve of generator $g$ at time $t$ in scenario $s$ (MW) & $\df_{e, t, s}$ & Thermal limit violation on branch $e$ at time $t$ in scenario $s$ (MW) \\
        \midrule
        \multicolumn{4}{l}{\textbf{Parameters}} \\
        \midrule
        $\dt$ & Duration of a time period (min) & $p_s$ & Probability of scenario $s$ \\
        $\rt$ & Duration of the ramping product (min) & $\pd_{i,t,s}$ & Power demand at bus $i$ at time $t$ in scenario $s$ (MW) \\
        $\sigma_{e,i}$ & PTDF coefficient for branch $e$ and bus $i$ & $\pgmin_{g,t,s}, \pgmax_{g,t,s}$ & Min and max output of generator $g$ at time $t$ in scenario $s$ (MW) \\
        $\rrup_{g,t}, \rrdn_{g,t}$ & Ramping-up and -down rate of generator $g$ at time $t$ (MW/min) & $\resReq_{t}$ & Reserve requirement (MW) \\
        $\rcupReq_{t}, \rcdnReq_{t}$ & Ramping-up and -down capability requirement at time $t$ (MW) & $\pfmin_{e,t}, \pfmax_{e,t}$ & Lower and upper thermal limits on branch $e$ at time $t$ (MW) \\
        $\costPg_{g, t}$ & Energy dispatch cost for generator $g$ at time $t$ (\$/MW) & $\costRes_{g, t}$ & Reserve cost for generator $g$ at time $t$ (\$/MW) \\
        $\penaltyRCUp, \penaltyRCDn$ & Ramping-up and -down capability shortage penalty price (\$/MW) & $\penaltyFlow_{e,t}$ & Thermal limit violation penalty price for branch $e$ at time $t$ (\$/MW) \\
        \bottomrule
    \end{tabularx}
\end{table}

\subsection{Variables}

The economic dispatch models comprise several key decision variables. 
The net power injection at bus $i$ for a given time $t$ and scenario $s$ is denoted by $\p_{i,t,s}$ and measured in megawatts (MW). 
Generator outputs are represented by $\pg_{g,t,s}$, indicating the energy dispatched by generator $g$ at time $t$ in scenario $s$.
For simplicity, different types of reserve capacities, such as regulating, operating, spinning and non-spinning reserves, are aggregated together and represented as $\res_{g,t,s}$ for each generator.
The model also accounts for the ramping capabilities of generators, which are specified by $\rcup_{g,t,s}$ for ramping up and $\rcdn_{g,t,s}$ for ramping down. Shortages in these ramping capabilities are captured by $\drcup_{t,s}$ and $\drcdn_{t,s}$, respectively. 
Lastly, $\df_{e,t,s}$ quantifies the thermal limit violations on branch $e$ at time $t$ in scenario $s$.

\subsection{Constraints}

\paragraph{Power balance constraints}

Equation \eqref{eq:stoch:nodal_power_balance} ensures nodal power balance by equating the net power at each node $i$ (the difference between generated power $\pg$ and demanded power $\pd$) for every time step $t$ and scenario $s$. Equation \eqref{eq:stoch:global_power_balance} enforces that the total net power across all nodes $\mathcal{N}$ sums to zero, maintaining the overall system balance for each time and scenario.
\begin{subequations}
    \setlength{\abovedisplayskip}{5pt}
    \setlength{\belowdisplayskip}{-10pt}
    \begin{align}
        \p_{i,t,s} & = \sum_{g \in \mathcal{G}_{i}} \pg_{g,t,s} - \pd_{i,t,s} , && \forall \, i,t,s \label{eq:stoch:nodal_power_balance}\\
        \sum_{i \in \mathcal{N}} \p_{i,t,s} & = 0, && \forall \, t,s \label{eq:stoch:global_power_balance}
    \end{align}
    \label{eq:stoch:power_balance}
\end{subequations}

\paragraph{Generation limits}
The generation constraints ensure that generator outputs remain within specified operational limits. Specifically, the sum of generated power $\pg_{g, t, s}$, reserved power $\res_{g, t, s}$, and ramp-up capability $\rcup_{g, t, s}$ must not exceed the maximum generation limit $\pgmax_{g,t,s}$ as shown in Equation \eqref{eq:stoch:eco_max}. Conversely, the difference between generated power, reserved power, and ramp-down capability $\rcdn_{g, t, s}$ must remain above the minimum generation limit $\pgmin_{g,t,s}$, ensuring operational reliability and flexibility as stated in Equation \eqref{eq:stoch:eco_min}.
\begin{subequations}
    \setlength{\abovedisplayskip}{0pt}
    \setlength{\belowdisplayskip}{0pt}
    \begin{align}
        \pg_{g, t, s} + \res_{g, t, s} + \rcup_{g, t, s} & \leq \pgmax_{g,t,s},  && \forall \, g,t,s \label{eq:stoch:eco_max}\\
        \pg_{g, t, s} - \res_{g, t, s} - \rcdn_{g, t, s} & \geq \pgmin_{g,t,s},  && \forall \, g,t,s \label{eq:stoch:eco_min}
    \end{align}
\end{subequations}

\paragraph{Generation ramp limits}
Ramping constraints for generators are critical to maintaining system stability and responsiveness. Equation \eqref{eq:stoch:gen:ramp_up} ensures that the power output of any generator $g$ at time $t$ in scenario $s$ must not exceed the previous time step's output incremented by the maximum ramp-up rate $\rrup_{g,t}$ scaled by the time duration $\dt$. Similarly, Equation \eqref{eq:stoch:gen:ramp_dn} ensures that the power output does not fall below the previous time step's output reduced by the maximum ramp-down rate $\rrdn_{g,t}$ scaled by the same time step $\dt$.
\begin{subequations}
    \setlength{\abovedisplayskip}{5pt}
    \setlength{\belowdisplayskip}{-15pt}
    \begin{align}                
        \pg_{g,t,s} & \leq \pg_{g,t-1,s} + \dt \, \rrup_{g,t} && \forall g, t, s \label{eq:stoch:gen:ramp_up}\\
        \pg_{g,t,s} & \geq \pg_{g,t-1,s} - \dt \, \rrdn_{g,t} && \forall g, t, s \label{eq:stoch:gen:ramp_dn}
    \end{align}
    \label{eq:stoch:gen:ramp}
\end{subequations}

\paragraph{Transmission Constraints}
Transmission constraints are essential for controlling the flow of electricity through transmission lines or transformers. They are formalized in Equations \eqref{eq:stoch:transmission_max} and \eqref{eq:stoch:transmission_min}, which ensure that the weighted sum of net power injections $\p_{i, t, s}$ at each bus $i$ does not exceed the maximum allowable flow $\pfmax_{e,t}$ and does not drop below the minimum allowable flow $\pfmin_{e,t}$. These limits are adjusted by $\df_{e, t, s}$ to account for operational deviations, helping to prevent overloads and maintain stable electricity distribution under varying operational scenarios.
\begin{subequations}
    \setlength{\abovedisplayskip}{5pt}
    \setlength{\belowdisplayskip}{-10pt}
    \begin{align}
        \sum_{i \in \mathcal{N}} \sigma_{e,i} \, \p_{i, t, s} & \leq \pfmax_{e,t} + \df_{e, t, s}, && \forall \, e,t,s \label{eq:stoch:transmission_max}\\
        \sum_{i \in \mathcal{N}} \sigma_{e,i} \, \p_{i, t, s} & \geq \pfmin_{e,t} - \df_{e, t, s}, && \forall \, e,t,s \label{eq:stoch:transmission_min}
    \end{align}
    \label{eq:stoch:transmission}
\end{subequations}

\paragraph{Non-negativity Constraints}
Equation \eqref{eq:stoch:non_negativity} ensures that reserves ($\res$), ramp capabilities ($\rcup$ and $\rcdn$), ramp capability shortages ($\drcup$ and $\drcdn$), and thermal limit deviations ($\df$) must all be nonnegative.
\setlength{\abovedisplayskip}{5pt}
\setlength{\belowdisplayskip}{5pt}
\begin{equation}
    \begin{aligned}
        \res, \rcup, \rcdn, \drcup , \drcdn , \df & \geq \mathbf{0}. \label{eq:stoch:non_negativity}
    \end{aligned}
\end{equation}

\paragraph{Ramping capability constraints}
\label{subsec:ramping_cap}
These constraints are also referred to as ``Flexiramp'' \citep{WangHobbs2013_Flexiramp}. Flexiramp is designed to ensure sufficient flexible generation capacity to handle volatile net loads effectively. These constraints are included only in the SCED+RP formulation and not in other formulations. Equations \eqref{eq:stoch:ramp_up_prod} and \eqref{eq:stoch:ramp_dn_prod} ensure that the total ramp-up and ramp-down capabilities of all generators meet the system's ramping requirements, $\rcupReq_{t}$ and $\rcdnReq_{t}$, respectively. These are modeled as soft constraints where $\drcup_{t, s}$ and $\drcdn_{t, s}$ represent the shortages in ramp-up and ramp-down capabilities. Additionally, Equations \eqref{eq:stoch:ramp_up_bound} and \eqref{eq:stoch:ramp_dn_bound} ensure that the ramping capabilities do not exceed the product of the specified ramp duration $\rt$ and the ramp rates. For instance, for a 30-minute ramping product, the ramping capabilities are limited by multiplying the 30-minute duration by the respective ramp rates in MW per minute.

\begin{subequations}
    \setlength{\abovedisplayskip}{-10pt}
    \setlength{\belowdisplayskip}{-10pt}
    \begin{align}
        \sum_{g \in \mathcal{G}} \left(\rcup_{g,t,s} \right) & \geq \rcupReq_{t} - \drcup_{t, s}, && \forall \, t,s\label{eq:stoch:ramp_up_prod} \\
        \sum_{g \in \mathcal{G}} \left(\rcdn_{g,t,s} \right) & \geq \rcdnReq_{t} - \drcdn_{t, s}, && \forall \, t,s\label{eq:stoch:ramp_dn_prod} \\
        \rcup_{g,t,s} & \leq \rt \, \rrup_{g,t}, && \forall \, g,t,s \label{eq:stoch:ramp_up_bound}\\
        \rcdn_{g,t,s} & \leq \rt \, \rrdn_{g,t}, && \forall \, g,t,s \label{eq:stoch:ramp_dn_bound}
    \end{align}
    \label{eq:stoch:ramp_cap}
\end{subequations}

\paragraph{Non-anticipatory Constraints}
In the stochastic formulation of our model, non-anticipatory constraints are crucial to ensure that decisions at the first time step remain consistent across all scenarios. Specifically, Equations \eqref{eq:stoch:non-anticipatory:pg}, \eqref{eq:stoch:non-anticipatory:res}, \eqref{eq:stoch:non-anticipatory:rcup}, and \eqref{eq:stoch:non-anticipatory:rcdn} require that the power generation ($\pg$), reserves ($\res$), and ramping capabilities ($\rcup$ and $\rcdn$) for any generator $g$ must be the same for all scenarios $\mathcal{S}$ in real-time operations ($t=1$).
\begin{subequations}
    \setlength{\abovedisplayskip}{0pt}
    \setlength{\belowdisplayskip}{-15pt}
    \begin{align}
        \pg_{g, 1, s_1} & = \pg_{g, 1, s_2}, && \forall \, g \in \mathcal{G}, s_1, s_2 \in \mathcal{S} \label{eq:stoch:non-anticipatory:pg}\\
        \res_{g, 1, s_1} & = \res_{g, 1, s_2}, && \forall \, g \in \mathcal{G}, s_1, s_2 \in \mathcal{S} \label{eq:stoch:non-anticipatory:res}\\
        \rcup_{g, 1, s_1} & = \rcup_{g, 1, s_2}, && \forall \, g \in \mathcal{G}, s_1, s_2 \in \mathcal{S} \label{eq:stoch:non-anticipatory:rcup}\\
        \rcdn_{g, 1, s_1} & = \rcdn_{g, 1, s_2}, && \forall \, g \in \mathcal{G}, s_1, s_2 \in \mathcal{S} \label{eq:stoch:non-anticipatory:rcdn}
    \end{align}
    \label{eq:stoch:non-anticipatory}
\end{subequations}

\subsection{Objective Function}
The objective function of the model is designed to minimize the total expected operational and penalty costs of the power system over all time periods $t$ in set $\mathcal{T}$, and across all scenarios $s$ in set $\mathcal{S}$. It is expressed in Equation \eqref{eq:stoch:obj} as the weighted sum of generation costs, reserve costs, and penalties associated with ramping capability shortages and thermal limit violations. 

\begin{equation}
    \sum_{t \in \mathcal{T}, s \in \mathcal{S}} p_{s} \dt \biggl(
    \sum_{g \in \mathcal{G}} \Bigl(
    \costPg_{g,t} \, \pg_{g,t,s} + \costRes_{g,t} \, \res_{g,t,s} \Bigr) 
    + \penaltyRCUp \, \drcup_{t,s}
    + \penaltyRCDn \, \drcdn_{t,s}
    + \sum_{e \in \mathcal{E}} \penaltyFlow_{e,t} \, \df_{e, t, s} \biggr)
    \label{eq:stoch:obj}
\end{equation}

\subsection{Deterministic Formulations}

\paragraph{Baseline RT-SCED}
\label{sec:formulation:SCED}

The baseline SCED formulation considered in the paper minimizes operational and penalty costs over a single time period and scenario.
It is formulated as the linear programming problem
\begin{align}
    \text{(SCED)} \quad
    \begin{array}{ll}
    \min \quad & \eqref{eq:stoch:obj} \\
    \text{s.t.} \quad & \eqref{eq:stoch:power_balance} - \eqref{eq:stoch:non_negativity} \quad \forall i \in \mathcal{N}, g \in \mathcal{G}, t \in \mathcal{T}, s \in \mathcal{S}
    \end{array}
    \label{eq:SCED}
\end{align}
where $\mathcal{T} = \{1\}$ and $\mathcal{S} = \{1\}$.

This formulation is the basis of, e.g., real-time markets in MISO \citep{Ma2009_MISO_RTSCED,MISO_BPM002_D}.
The objective \eqref{eq:stoch:obj} computes the overall cost of operation and penalties.
The constraints capture global power balance, generator and reserve min/max limits, reserve requirements, and transmission constraints. 
Additionally, the formulation includes ramping up and down requirements to ensure that changes in dispatch can be accommodated by the ramp-up rate.

\subsubsection{RT-SCED with Flexiramp}

Following the ramping products design outlined by MISO \citep{Chen2023_MISORPDesign}, the SCED with ramping products (SCED+RP) formulation incorporates 10-minute ramping-up and ramping-down capability products, along with 30-minute short-term ramping products. The ramping requirements are calculated based on changes and uncertainties in net load. Model \eqref{eq:SCEDRP} illustrates the corresponding formulation. Compared to Model \eqref{eq:SCED}, this formulation includes ramping capability constraints described in \ref{subsec:ramping_cap}, also known as ramping products.
\begin{align}
    \text{(SCED+RP)} \quad
    \begin{array}{ll}
    \min \quad & \eqref{eq:stoch:obj} \\
    \text{s.t.} \quad & \eqref{eq:stoch:power_balance} - \eqref{eq:stoch:ramp_cap}, \quad \forall i \in \mathcal{N}, g \in \mathcal{G}, t \in \mathcal{T}, s \in \mathcal{S}
    \end{array}
    \label{eq:SCEDRP}
\end{align}
where $\mathcal{T} = \{1\}, \mathcal{S} = \{1\}$.

\subsubsection{Look-Ahead Dispatch}
\label{sec:formulation:LAD}

Look-Ahead Dispatch (LAD) extends single-period SCED to a
multi-period economic dispatch formulation over a horizon of length
$T$.
The corresponding compact formulation is
\begin{align}
    \text{(LAD)} \quad
    \begin{array}{ll}
    \min \quad & \eqref{eq:stoch:obj} \\
    \text{s.t.} \quad & \eqref{eq:stoch:power_balance} - \eqref{eq:stoch:non_negativity} \quad \forall i \in \mathcal{N}, g \in \mathcal{G}, t \in \mathcal{T}, s \in \mathcal{S}
    \end{array}
    \label{eq:LAD}
\end{align}
where $\mathcal{T} = \{1, ..., T\}$ and $\mathcal{S} = \{1\}$.

The problem takes as inputs deterministic forecasts for load and
renewable generation, and outputs the energy dispatch and reserves for
each generator and each period.  Transmission constraints and reserve
requirements are also enforced for each period. Decisions at time $t$
are linked to decisions at time $t+1$ via ramping constraints for the
generators (constraints \eqref{eq:stoch:gen:ramp_up} and
\eqref{eq:stoch:gen:ramp_dn}). If the ramping constraints are not binding,
LAD decomposes into $T$ independent SCED problems, and the
solution at $t \, {=} \, 1$ is equivalent to the single-period SCED solution.
    
\subsection{Stochastic Formulation}
\label{sec:formulation:SLAD}

Stochastic look-ahead dispatch (SLAD) further extends LAD by considering uncertainty in future electricity demand and renewable generation, i.e., SLAD considers multiple time steps and multiple scenarios.
Namely, the formulation considered in this paper is a two-stage formulation, where the first-stage decisions are the decision variables at $t \, {=} \, 1$ and the second-stage decisions are decision variables at times $t\, {=} \, 2, ..., T$.
This reflects the fact that only dispatch decisions for $t=1$ are communicated to generators and implemented; subsequent decisions ($t=2, ..., T$) are only advisory.
The resulting extended SLAD formulation reads
\begin{align}
    \text{(SLAD)} \quad
    \begin{array}{ll}
    \min \quad & \eqref{eq:stoch:obj} \\
    \text{s.t.} \quad & \eqref{eq:stoch:power_balance} - \eqref{eq:stoch:non_negativity}, \eqref{eq:stoch:non-anticipatory} \quad \forall i \in \mathcal{N}, g \in \mathcal{G}, t \in \mathcal{T}, s \in \mathcal{S}
    \end{array}
    \label{eq:SLAD}
\end{align}
where $\mathcal{T}  \, {=} \, \{1, ..., T\}$ and $\mathcal{S} \, {=} \, \{1, ..., S\}$.
Naturally, if a single scenario is considered, i.e., $| \mathcal{S}| \, {=} \,1$, SLAD is equivalent to LAD.
Table \ref{tab:abstract_formulations} summarizes the ED formulations considered in the paper.

\begin{table}[!t]
    \centering
    \caption{The Economic Dispatch Formulations.}
    \label{tab:abstract_formulations}
    \begin{tabular}{cccccc}
        \toprule
            Problem & Time steps ($|\mathcal{T}|$) & Scenarios ($|\mathcal{S}|$) & FRPs$^{*}$\\
        \midrule
            SCED    & \phantom{$>$}1  & \phantom{$>$}1     & \xmark \\
            SCED+RP & \phantom{$>$}1  & \phantom{$>$}1     & \cmark \\
            LAD     & $>$1          & \phantom{$>$}1     & \xmark \\
            SLAD    & $>$1          & $>$1             & \xmark \\
        \bottomrule
    \end{tabular}
    \par
    \footnotesize{$^{*}$FRPs: Flexible Ramp Products.}
\end{table}

\emph{Non-anticipatory} constraints \eqref{eq:stoch:non-anticipatory} ensure that the first-stage decisions $\mathbf{x}_{1,s}$ are the same across all scenarios.
Problem \eqref{eq:SLAD} is a (large) linear problem, with $\mathcal{O}(T {\times} S)$ variables and constraints, which is
tractable only for small to medium-sized instances.  
Large instances require specialized algorithms such as Benders decomposition.

\section{Accelerated Benders Decomposition}
\label{sec:method}

The extensive formulation \eqref{eq:SLAD} is not tractable for the
large-scale instances considered in this work, especially for
real-time operations. Instead, specialized algorithms are needed, such
as Benders Decomposition (BD) \citep{Benders1962_Decomposition} and the
Progressive Hedging Algorithm (PHA) \citep{Rockafellar1991_ProgressiveHedging}.
In preliminary experiments, Benders decomposition was found to be faster and exhibit better convergence than PHA, hence, the paper focuses on BD.
This section presents the Benders decomposition (BD) algorithm and acceleration techniques to improve its convergence for large SLAD instances.

For ease of reading, the proposed accelerated Benders decomposition is stated for a general two-stage stochastic program of the form
\begin{subequations}
    \setlength{\abovedisplayskip}{5pt}
    \setlength{\belowdisplayskip}{5pt}
    \label{eq:BD:TSSP}
    \begin{align}
        \min_{\x, \y} \quad 
            & c^{\top} \x + \sum_{s \in \mathcal{S}} q_{s}^{\top} \y_{s} \label{eq:BD:TSSP:obj}\\
        \text{s.t.} \quad
            & A \x \geq b, \label{eq:BD:TSSP:con:first_stage}\\
            & H_{s} \x + W_{s} \y_{s} \geq h_{s}, \quad \forall s \in \mathcal{S} \label{eq:BD:TSSP:con:second_stage}
    \end{align}
\end{subequations}
where $x$ and $y$ denote the first-stage and second-stage variables, respectively.
The objective \eqref{eq:BD:TSSP:obj} minimizes the sum of first-stage and expected second-stage costs; for compactness, the probability associated with scenario $s$ is captured in the definition of $q_{s}$.
Constraints \eqref{eq:BD:TSSP:con:first_stage} and \eqref{eq:BD:TSSP:con:second_stage} capture first-stage and second-stage constraints, respectively.
Recall that, in the context of SLAD, first-stage variables $\x$ comprise all decision variables for the first time step $t=1$, and second-stage variables $y$ comprise all decision variables for time steps $t=2, ..., T$.

\subsection{Benders Decomposition}
\label{sec:method:BD}

    The Benders master problem considers only the first-stage variables
    and is of the form
    \begin{subequations}
        \setlength{\abovedisplayskip}{5pt}
        \setlength{\belowdisplayskip}{5pt}
        \label{eq:BD:master}
        \begin{align}
            \min_{\x, \theta} \quad 
                & c^{\top} \x + \sum_{s \in \mathcal{S}} \theta_{s} \label{eq:BD:master:obj}\\
            \text{s.t.} \quad 
                & A \x \geq b, \label{eq:BD:master:first_stage}\\
                & \theta_{s} \geq Q_{s}(\x), && \forall s \in \mathcal{S}\label{eq:BD:master:second_stage}
        \end{align}
    \end{subequations}
    where $Q_{s}(\x)$ denotes the second-stage cost of $\x$ under scenario $s$, defined as
    \begin{subequations}
        \setlength{\abovedisplayskip}{5pt}
        \setlength{\belowdisplayskip}{5pt}
        \label{eq:BD:SP}
        \begin{align}
            Q_{s}(\x) =
            \min_{\y_{s}} \quad
                & q_{s}^{\top} \y_{s}\\
            \text{s.t.} \quad
                & W_{s} \y_{s} \geq h_{s} - H_{s} \x. \label{eq:BD:SP:con}
        \end{align}
    \end{subequations}
    In all economic dispatch formulations considered in the paper, system-wide constraints such as power balance, reserve requirements and transmission constraints, are treated as soft constraints.
    Therefore, SLAD has relatively complete recourse, i.e., for any first-stage solution $\bar{\x}$ that satisfies constraints the sub-problem \eqref{eq:BD:SP} is feasible 
    
    The second-stage cost function $Q_{s}$ is the value function of a linear programming problem, hence, it is a convex, piece-wise linear function.
    However, in general, $Q_{s}$ contains exponentially many pieces, making its direct representation intractable.
    Instead, Benders decomposition uses a cutting-plane approach to iteratively refine a piece-wise linear under-approximation of $Q_{s}$ until convergence.
    Denoting by $\lambda_{s}$ the dual variable associated with constraint \eqref{eq:BD:SP:con}, each Benders cut has the form
    \begin{align}
        \label{eq:BD:cut}
        \theta_{s} \geq \lambda_{s}^{\top} (h_{s} - H_{s} \x).
    \end{align}
    The cuts highlighted in \eqref{eq:BD:cut} are also referred to as \emph{Benders optimality cuts}.
    Because SLAD has relatively complete recourse, it is sufficient to consider optimality cuts.
    The reader is referred to \citep{Benders1962_Decomposition} for the general treatment of so-called feasibility cuts, which are needed when subproblems can be infeasible.
    In what follows, using a slight abuse of notation, each BD cut is identified with the corresponding dual solution $\lambda_{s} \in \Lambda_{s}$ where, $\forall s \in \mathcal{S}$, $\Lambda_{s}$ denotes the dual-feasible set
    \begin{align}
        \label{eq:BD:dual_feasible}
        \Lambda_{s} = \left\{
            \lambda_{s} \geq 0
        \ \middle| \
            W_{s}^{\top}\lambda_{s} = q_{s}
        \right\}.
    \end{align}

    The Benders master problem is initialized with a small number of cuts $\bar{\Lambda} \subseteq \Lambda = (\Lambda_{1}, ..., \Lambda_{S})$.
    The corresponding relaxed master problem (RMP) reads
    \begin{subequations}
        \setlength{\abovedisplayskip}{5pt}
        \setlength{\belowdisplayskip}{5pt}
        \label{eq:BD:RMP}
        \begin{align}
            RMP(\bar{\Lambda}) \quad \quad \min_{\x, \theta} \quad 
                & c^{\top} \x + \sum_{s \in \mathcal{S}} \theta_{s} \label{eq:BD:RMP:obj}\\
            \text{s.t.} \quad 
                & A \x \geq b, \label{eq:BD:RMP:first_stage}\\
                & \theta_{s} \geq \lambda_{s}^{\top}(h_{s} - H_{s} \x),
                    && \forall s \in \mathcal{S}, \forall \lambda_{s} \in \bar{\Lambda}_{s}.
                    \label{eq:BD:RMP:cuts}
        \end{align}
    \end{subequations}
    At each iteration, new cuts are identified by solving the Benders sub-problem \eqref{eq:BD:SP} for each scenario, and retrieving the corresponding optimal dual solution.
    Violated cuts are added to the master problem, and the algorithm terminates with an optimal solution when no violated cut exists.
    Although BD is guaranteed to converge in a finite number of iterations, textbook implementations often suffer from slow convergence.
    Hence, the paper implements several acceleration strategies described below.

\subsection{Acceleration Strategies for the Subproblems}
\label{sec:BD:SP}

    Two mechanisms are implemented to speed up the resolution of the BD subproblems.
    First, using shared or distributed memory architectures, the subproblems are solved in parallel, each subproblem being handled by an independent worker process.
    At each iteration, the candidate master solution is sent to each worker process, which returns the solution of the corresponding subproblem.
    This strategy yields significant speedups, especially since, for the case at hand, optimizing the subproblems is substantially more time-consuming than optimizing the master problem.

    Second, each subproblem contains a large number of transmission constraints of the form \eqref{eq:stoch:transmission}, which are modeled via a Power Transfer Distribution Factor (PTDF) formulation.
    Because only a small number of transmission lines are congested in practice, PTDF constraints are handled in a lazy fashion.
    Namely, each subproblem is initialized without any transmission constraints, and solved in an iterative fashion.
    Thereby, at each iteration, the current solution is checked for violations, and only violated transmission constraints are added to the problem.
    This process is repeated until all transmission constraints are satisfied.
    By leveraging the warm-start capabilities of the simplex algorithm, this strategy significantly speeds up the resolution of individual subproblems.
    Importantly, it does not impact the validity nor quality of Benders cuts.
    For the same reason, transmission constraints are also handled in lazy fashion when solving the Benders master problem.

\subsection{In-out Separation}
\label{sec:BD:in-out}
    
    To alleviate the slow convergence and tailing off effect of a vanilla BD, the paper implements an \textit{in-out separation} technique \citep{fischetti2016benders, fischetti2010inout}.
    Thereby, at each iteration, BD cuts are generated by seeding the subproblems with a convex combination of the current master solution $\bar{\x}$ and a core point $\hat{\x}$, which resides in the interior of the feasible region.
    Namely, the subproblems receive $\tilde{\x} = \alpha \bar{\x} + (1-\alpha) \hat{\x}$, where $\alpha \in [0, 1]$ is typically chosen by performing a line search over the segment $[0, 1]$.
    This strategy is known to generate stronger cuts \citep{bonami_2020}, thereby improving convergence, especially in later iterations.

    In practice, one should be mindful of the following facts regarding the choice of $\alpha$.
    On the one hand, a full line search comes at the cost of having to solve each subproblem multiple times per BD iteration, which may be computationally burdensome.
    On the other hand, using a fixed $\alpha < 1$ throughout the algorithm does not guarantee convergence.

    To balance these issues, the paper uses the following strategy.
    At each BD iteration, subproblems are first seeded with $\tilde{\mathbf{x}} \, {=} \, \bar{\alpha} \bar{\mathbf{x}} \, {+} \, (1 \, {-} \, \bar{\alpha}) \hat{\mathbf{x}}$, where $0 \, {<} \, \bar{\alpha} \, {<} \, 1$ is a user-specified parameter.
    If a violated cut is found, it is added to the master problem and the algorithm continues to next BD iteration.
    If no violated cut is found, the subproblems are re-seeded with the current master solution $\bar{\mathbf{x}}$, which is equivalent to setting $\alpha \, {=} \, 1$.
    This scheme guarantees finite convergence while ensuring that subproblems are solved at most twice at each iteration.
    Indeed, if no violated cut is found after the second step, then the current master solution $\bar{x}$ is optimal and the algorithm terminates.

    Another important aspect of the in-out strategy is the choice of core point $\hat{\x}$.
    In the paper's implementation, $\hat{\mathbf{x}}$ is initialized as the solution of the initial master problem.
    Then, at each iteration, $\hat{\x}$ is set to the last solution that was sent to the subproblems.
    Thereby, $\hat{\x}$ is set to $\tilde{\x}$ if the first round of separation yielded a violated cut, and to $\bar{\x}$ otherwise.
    This aggressive strategy was found to be effective at preventing tailing off.
    Finally, note that, although this scheme does not guarantee that $\hat{\x}$ lies in the interior of the feasible region, it does not affect the validity of the algorithm.

\subsection{Accelerated Benders Decomposition Algorithm}
\label{sec:BD:algorithm}

    \begin{algorithm}[!t]
        \caption{Accelerated Benders Decomposition
        }
        \label{alg:benders}
        \begin{algorithmic}[1]
            \REQUIRE {Tolerance $\epsilon \geq 0$; in-out parameter $\bar{\alpha} \in (0, 1)$; initial set of cuts $\bar{\Lambda} = (\bar{\Lambda}_{1}, ..., \bar{\Lambda}_{S})$}
            \STATE $(LB, UB) \gets (-\infty, + \infty)$
            \STATE Solve $\RMP{\bar{\Lambda}}$; retrieve optimal solution $\bar{\x}$ \label{alg:BD:masater}
            \STATE $\hat{\x} \gets \bar{\x}$ \label{alg:BD:in_out:init}
            \WHILE{$(UB - LB) \geq \epsilon$} {
                \STATE Solve $\RMP{\bar{\Lambda}}$; retrieve optimal solution $(\bar{\x}, \bar{\theta})$ \label{alg:BD:solve_rmp}
                \STATE $LB \gets \max(LB, c^{\top} \bar{\x} + \sum_{s} \bar{\theta}_{s}) $ \label{alg:BD:lb}
                \STATE $(\alpha, n^{c}) \gets (\bar{\alpha}, 0)$
                \WHILE{$n^{c} = 0$} {
                    \STATE $\tilde{\x} \gets \alpha \bar{\x} + (1-\alpha)\hat{\x}$ \label{alg:BD:cut_loop:cvx}
                    \FOR {$s = 1,2,\dots, S$ \label{alg:BD:cut_loop:begin}}{
                        \STATE {Solve subproblem \eqref{eq:BD:SP} to obtain $Q_{s}(\tilde{\x})$} \label{alg:BD:cut_loop:subproblem} and dual solution $\lambda_{s}$ \label{alg:BD:cut_loop:sp_dual}
                        \IF {$\bar{\theta}_{s} < \lambda_{s}^{\top}(h_{s} - H_{s} \bar{\x})$}{
                            \STATE $\bar{\Lambda}_{s} \gets \bar{\Lambda}_{s} \cup \{ \lambda_{s} \}$ \label{alg:BD:cut_loop:add_cut}
                            \STATE $n^{c} \gets n^{c} + 1$
                        } \ENDIF
                    } \ENDFOR \label{alg:BD:cut_loop:end}
                    \IF {($n^{c} = 0)$ and $(\alpha < 1)$}{
                        \STATE $\alpha \gets 1$ \label{alg:BD:in_out:reinit}
                    } \ELSIF {($n^{c} = 0)$ and $(\alpha = 1)$}{
                        \STATE break
                    } \ENDIF
                } \ENDWHILE
                \STATE $UB \gets \min \left(
                    UB,
                    c^{\top} \tilde{\x} + \sum_{s} Q_{s}(\tilde{\x})
                    \right)
                $ \label{alg:BD:ub}
                \STATE $\hat{\x} \gets \tilde{\x}$
            } \ENDWHILE
        \end{algorithmic}
    \end{algorithm}

    The proposed Accelerated Benders Decomposition (ABD) algorithm, incorporating parallel solving of sub-problems, lazy thermal constraints, and in-out separation, is detailed in Algorithm \ref{alg:benders}.

    To prevent the relaxed master problem from being unbounded, it is initialized with trivial valid cuts of the form $\theta_{s} \, {\geq} \, M, \forall s$, where $M$ is a large negative constant.
    The core point $\hat{\mathbf{x}}$ is initialized at line \ref{alg:BD:in_out:init} using the initial RMP solution.
    Each ABD iteration solves the RMP to obtain $(\bar{\mathbf{x}}, \bar{\theta})$ (line \ref{alg:BD:solve_rmp}), updates the lower bound (line \ref{alg:BD:lb}), and enters an in-out separation loop.
    This loop generates Benders cuts by seeding the subproblems with $\tilde{\mathbf{x}} = \alpha \bar{\mathbf{x}} + (1 - \alpha) \hat{\mathbf{x}}$ (line \ref{alg:BD:cut_loop:cvx}) using a convex weight $\alpha$.
    Subproblems are solved for each scenario using $\tilde{\mathbf{x}}$ (lines \ref{alg:BD:cut_loop:begin}--\ref{alg:BD:cut_loop:end}), and violated cuts are added to the RMP (line \ref{alg:BD:cut_loop:add_cut}).
    If no violated cuts are found during the first separation round (i.e. $n^{c} \, {=} \, 0$ and $\alpha \, {=} \, 1$), $\alpha$ is set to 1 at line \ref{alg:BD:in_out:reinit}, effectively moving the separation point $\tilde{\x}$ to the current master solution $\bar{\x}$ and triggering a re-solve of the sub-problems.
    If no violated cuts are found after setting $\alpha$ to $1$, the separation loop is interrupted (line 20); note that this indicates that the current solution is optimal.
    This scheme ensures that subproblems are solved at most twice per iteration while retaining the validity of the overall algorithm.
    Finally, the upper bound is updated based on subproblem solutions (line \ref{alg:BD:ub}), and the core point $\hat{\x}$ is updated.
    This iterative process is repeated until the optimality gap $(UB - LB)$ reaches the prescribed tolerance $\epsilon$.

\section{Illustrative Example}
\label{sec:example}

\colorlet{SCED}{Purple}
\colorlet{SCEDRP}{red}
\colorlet{LAD}{orange}
\colorlet{SLAD}{blue}

This section illustrates the behavior of each Economic Dispatch (ED) formulation
on a small system over two time periods.  The system at hand comprises
two 20MW generators $G_{1}$ and $G_{2}$ and a single load; for simplicity,
no reserve products nor transmission lines are considered.  Table
\ref{tab:illustrative_example:gen_info} reports the operating
parameters of both generators.  Generator $G_{1}$ has lower cost and higher ramping rate; generator $G_{2}$ has higher cost, and lower ramping rate.

\begin{table}[!t]
    \centering
    \caption{Generator operating parameters}
    \label{tab:illustrative_example:gen_info}
    \begin{tabular}{crrr}
        \toprule
        Unit 
            & \multicolumn{1}{c}{Capacity} 
            & \multicolumn{1}{c}{Cost}
            & \multicolumn{1}{c}{Ramp rate}\\ 
        \midrule
        $G_1$ & 20 MW   & \$10/MW  & 4 MW/min \\
        $G_2$ & 20 MW   & \$20/MW  & 2 MW/min \\
        \bottomrule
    \end{tabular}
\end{table}

The look-ahead horizon in LAD and SLAD formulations is $T\,{=}\,2$, and the number of scenarios in SLAD is $S\,{=}\,2$.
Table \ref{tab:illustrative_example:demand_and_solution} reports the demand forecast used in each formulation for $t,{=},1, 2$. 
At each time step, SCED and SCED+RP only consider the current known demand, which is the actual demand, whereas LAD and SLAD require future demand information and therefore rely on forecasts.
SCED+RP additionally considers a 22MW requirement for a 5-minute flexible ramping product at $t\,{=}\,1$.
The deterministic forecast in LAD underestimates future demand.
Finally, note that, at each time, the average of the two scenarios considered in SLAD is equal to the deterministic forecast used by LAD.

\begin{table}[!t]
    \centering
    \caption{Forecast demand and solutions for each ED formulation}
    \label{tab:illustrative_example:demand_and_solution}
    \begin{tabular}{l cccc cccc}
        \toprule  
            & \multicolumn{4}{c}{$t = 1$}
            & \multicolumn{4}{c}{$t = 2$}\\
        \cmidrule(lr){2-5} \cmidrule(lr){6-9}
            Problem 
            & Forecast & $\pg$ & $\delta$ & Cost
            & Forecast & $\pg$ & $\delta$ & Cost$^{\dagger}$\\
        \midrule
        SCED 
            & 10 & (10, 0) & 0 & 100 
            & 35 & (20, 10) & 5 & 5400\\
        SCED+RP 
            & 10 & (8, 2) & 0 & 120 
            & 35 & (20, 12) & 3 & 3440\\
        LAD 
            & (10, 33) & (7, 3) & 0 & 130
            & (35, 29) & (20, 13) & 2 & 2460\\
        SLAD$^{*}$ 
            & $\begin{pmatrix} 10 & 29\\ 10 & 37\end{pmatrix}$ 
            & (3, 7) & 0 & 170
            & $\begin{pmatrix} 35 & 27\\ 35 & 31\end{pmatrix}$ 
            & (20, 15) & 0 & 500\\
        \bottomrule
    \end{tabular}\\
    $^{*}$each row corresponds to one scenario.
    \footnotesize{$^{\dagger}$include \$1000 per MW of shortage.}
\end{table}

Table \ref{tab:illustrative_example:demand_and_solution} reports, for each formulation and each
time $t\,{=}\,1,2$: the active power dispatch of both generators ($\pg \,
{\in} \, \mathbb{R}^{2}$, in MW), the corresponding energy shortage
($\delta$, in MW), and the operating costs (in \$).  The latter
comprises energy dispatch costs at time $t$, and includes a penalty of
\$1000 per MW of energy shortage.  Note that, for multi-period
formulations, only the first-stage decision is considered when
computing actual operating costs.

The results of Table \ref{tab:illustrative_example:demand_and_solution} show that SCED achieves the lowest costs at $t \, {=} \, 1$, but the highest total costs because of energy shortages at $t \, {=} \, 2$, which result in large penalty costs.
Similarly, SCED+RP and LAD also achieve lower costs at $t \, {=} \, 1$, but experiences energy shortages at $t \, {=} \, 2$, albeit less than SCED.
On the other hand, SLAD, despite incurring a higher cost at $t \, {=} \, 1$, avoids energy shortages at $t \, {=} \, 2$ and therefore achieves lower costs over the horizon.
This illustrates the fact that the short-term focus of SCED results in sub-optimal decisions in the long term.

\tikzset{
    scedstyle/.style={rectangle, draw=SCED, fill=SCED, inner sep=0pt, minimum size=3pt},
    scedrpstyle/.style={diamond, draw=SCEDRP, fill=SCEDRP, inner sep=0pt, minimum size=4pt},
    ladstyle/.style={regular polygon, regular polygon sides=3, draw=LAD, fill=LAD, inner sep=0pt, minimum size=4pt},
    sladstyle/.style={circle, draw=SLAD, fill=SLAD, inner sep=0pt, minimum size=3pt}
}

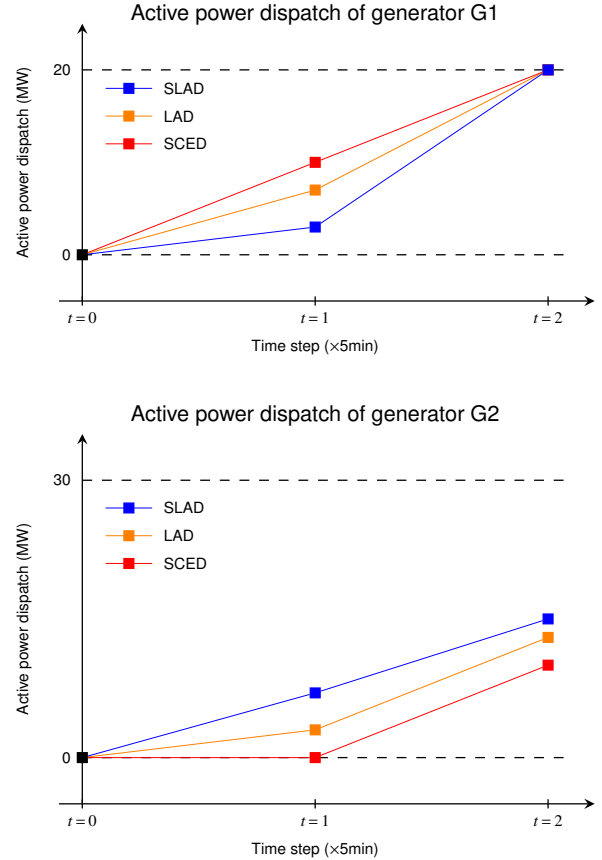
\begin{figure}[!t]
    \centering
    \resizebox{0.48\textwidth}{!}{
    \begin{tikzpicture}[x=2.5cm,y=0.12cm,font=\sffamily]
        \draw[black,-stealth] (-0.1, -5) -- (2.2, -5);
        \node [] at (1, -10) {\tiny{Time step ($\times$5min)}};
        \draw[black,-stealth] (0.0, -5) -- (0.0, 25);
        \node [rotate=90] at (-0.25, 10) {\tiny{Active power dispatch (MW)}};
        
        \node [] at (1.05, 26) {\scriptsize{Active power dispatch of generator G1}};

        \draw[black] (0.0, -5.5) -- (0.0, -4.5);
        \draw[black] (1.0, -5.5) -- (1.0, -4.5);
        \draw[black] (2.0, -5.5) -- (2.0, -4.5);
        \node[below] at (0.0, -5) {\tiny{$t \,{=} \, 0$}};
        \node[below] at (1.0, -5) {\tiny{$t \,{=} \, 1$}};
        \node[below] at (2.0, -5) {\tiny{$t \,{=} \, 2$}};
        
        \draw[black, dashed] (0,00) -- (2.1, 00);  
        \node[left] at (0, 0) {\tiny{0}};
        \draw[black, dashed] (0,20) -- (2.1, 20);  
        \node[left] at (0, 20) {\tiny{20}};

        \node[draw=black, fill=black, inner sep=0pt,minimum size=3pt] (pg1_0) at (0, 0) {};  

        \node[scedstyle] (pg1_1_sced) at (1, 10) {};  
        \node[scedrpstyle] (pg1_1_scedrp) at (1, 8) {};  
        \node[ladstyle] (pg1_1_lad)  at (1,  7) {};  
        \node[sladstyle] (pg1_1_slad) at (1,  3) {};  

        \draw[SCED, line width=0.25pt] (pg1_0) -- (pg1_1_sced);  
        \draw[SCEDRP, line width=0.25pt] (pg1_0) -- (pg1_1_scedrp);  
        \draw[LAD, line width=0.25pt] (pg1_0) -- (pg1_1_lad);  
        \draw[SLAD, line width=0.25pt] (pg1_0) -- (pg1_1_slad);  

        \node[scedstyle] (pg1_2_sced) at (2, 20) {};  
        \node[scedrpstyle] (pg1_2_scedrp) at (2, 20) {};  
        \node[ladstyle] (pg1_2_lad) at (2, 20) {};  
        \node[sladstyle] (pg1_2_slad) at (2, 20) {};  
        
        \draw[SCED, line width=0.25pt] (pg1_1_sced) -- (pg1_2_sced);  
        \draw[SCEDRP, line width=0.25pt] (pg1_1_scedrp) -- (pg1_2_scedrp);  
        \draw[LAD, line width=0.25pt] (pg1_1_lad) -- (pg1_2_lad);  
        \draw[SLAD, line width=0.25pt] (pg1_1_slad) -- (pg1_2_slad);  

        \node[scedstyle] at (0.2, 18) {};
        \draw[SCED, line width=0.25pt] (0.1, 18) -- (0.3, 18);
        \node[right] at (0.3, 18) {\tiny{SCED}};
        
        \node[scedrpstyle] at (0.2, 15) {};
        \draw[SCEDRP, line width=0.25pt] (0.1, 15) -- (0.3, 15);
        \node[right] at (0.3, 15) {\tiny{SCED+RP}};
        
        \node[ladstyle] at (0.2, 12) {};
        \draw[LAD, line width=0.25pt] (0.1, 12) -- (0.3, 12);
        \node[right] at (0.3, 12) {\tiny{LAD}};
        
        \node[sladstyle] at (0.2, 9) {};
        \draw[SLAD, line width=0.25pt] (0.1, 9) -- (0.3, 9);
        \node[right] at (0.3, 9) {\tiny{SLAD}};
    \end{tikzpicture}
    }
    \hfill
    \resizebox{0.48\textwidth}{!}{
    \begin{tikzpicture}[x=2.5cm,y=0.12cm,font=\sffamily]
        \draw[black,-stealth] (-0.1, -5) -- (2.2, -5);
        \node [] at (1, -10) {\tiny{Time step ($\times$5min)}};
        \draw[black,-stealth] (0.0, -5) -- (0.0, 25);
        \node [rotate=90] at (-0.25, 10) {\tiny{Active power dispatch (MW)}};
        
        \node [] at (1.05, 26) {\scriptsize{Active power dispatch of generator G2}};

        \draw[black] (0.0, -5.5) -- (0.0, -4.5);
        \draw[black] (1.0, -5.5) -- (1.0, -4.5);
        \draw[black] (2.0, -5.5) -- (2.0, -4.5);
        \node[below] at (0.0, -5) {\tiny{$t \,{=} \, 0$}};
        \node[below] at (1.0, -5) {\tiny{$t \,{=} \, 1$}};
        \node[below] at (2.0, -5) {\tiny{$t \,{=} \, 2$}};
        
        \draw[black, dashed] (0,00) -- (2.1, 00);  
        \node[left] at (0, 0) {\tiny{0}};
        \draw[black, dashed] (0,20) -- (2.1, 20);  
        \node[left] at (0, 20) {\tiny{20}};

        \node[draw=black, fill=black, inner sep=0pt,minimum size=3pt] (pg2_0) at (0, 0) {};  

        \node[scedstyle] (pg2_1_sced) at (1, 0) {};  
        \node[scedrpstyle] (pg2_1_scedrp) at (1, 2) {};  
        \node[ladstyle] (pg2_1_lad) at (1, 3) {};  
        \node[sladstyle] (pg2_1_slad) at (1, 7) {};  
        \draw[SCED, line width=0.25pt] (pg2_0) -- (pg2_1_sced);  
        \draw[SCEDRP, line width=0.25pt] (pg2_0) -- (pg2_1_scedrp);  
        \draw[LAD, line width=0.25pt] (pg2_0) -- (pg2_1_lad);  
        \draw[SLAD, line width=0.25pt] (pg2_0) -- (pg2_1_slad);  
        
        \node[scedstyle] (pg2_2_sced) at (2, 10) {};  
        \node[scedrpstyle] (pg2_2_scedrp) at (2, 12) {};  
        \node[ladstyle] (pg2_2_lad) at (2, 13) {};  
        \node[sladstyle] (pg2_2_slad) at (2, 15) {};  
        \draw[SCED, line width=0.25pt] (pg2_1_sced) -- (pg2_2_sced);  
        \draw[SCEDRP, line width=0.25pt] (pg2_1_scedrp) -- (pg2_2_scedrp);  
        \draw[LAD, line width=0.25pt] (pg2_1_lad) -- (pg2_2_lad);  
        \draw[SLAD, line width=0.25pt] (pg2_1_slad) -- (pg2_2_slad);  

        \node[scedstyle] at (0.2, 18) {};
        \draw[SCED, line width=0.25pt] (0.1, 18) -- (0.3, 18);
        \node[right] at (0.3, 18) {\tiny{SCED}};
        
        \node[scedrpstyle] at (0.2, 15) {};
        \draw[SCEDRP, line width=0.25pt] (0.1, 15) -- (0.3, 15);
        \node[right] at (0.3, 15) {\tiny{SCED+RP}};
        
        \node[ladstyle] at (0.2, 12) {};
        \draw[LAD, line width=0.25pt] (0.1, 12) -- (0.3, 12);
        \node[right] at (0.3, 12) {\tiny{LAD}};
        
        \node[sladstyle] at (0.2, 9) {};
        \draw[SLAD, line width=0.25pt] (0.1, 9) -- (0.3, 9);
        \node[right] at (0.3, 9) {\tiny{SLAD}};
        
    \end{tikzpicture}
    }
    \caption{Active power dispatch of generator G1 (top) and G2 (bottom) time  $t \, {=} \, 0, 1, 2$, under each economic dispatch formulation (\textcolor{SCED}{\textbf{SCED}}/\textcolor{SCEDRP}{\textbf{SCED+RP}}/\textcolor{LAD}{\textbf{LAD}}/\textcolor{SLAD}{\textbf{SLAD}}). The black square is the initial dispatch ($t \, {=} \, 0$).}
    \label{fig:example:dispatch}
\end{figure}

\begin{figure}[!t]
    \centering
    \resizebox{0.6\columnwidth}{!}{
    \begin{tikzpicture}[x=0.15cm,y=0.15cm,font=\sffamily]
        \draw[black,-stealth] (0, 0.0) -- (22, 0.0);
        \node[] at (10, -3) {\tiny{Active power dispatch of G1}};
        \draw[black,-stealth] (0.0, 0) -- (0.0, 22);
        \node[rotate=90] at (-3, 10) {\tiny{Active power dispatch of G2}};
        
        \draw[black] (0,0) -- (0, 20) -- (20, 20) -- (20, 0) -- cycle;
        
        \node[left] at (0, 0) {\tiny{$t \,{=}\, 0$}};
        \draw[black,line width=0.5pt, dash pattern=on 2pt off 1pt] (0, 10) -- (10, 0);  
        \node[right] at (0, 10) {\tiny{$t \,{=}\, 1$}};
        \draw[black,line width=0.5pt, dash pattern=on 2pt off 1pt] (15, 20) -- (20, 15);  
        \node[left] at (17, 18) {\tiny{$t \,{=}\, 2$}};

        \node[draw=black, fill=black, inner sep=0pt,minimum size=3pt] (pg_0) at (0, 0) {};
        
        \node[scedstyle] (pg_1_sced) at (10, 0) {};  
        \node[scedrpstyle] (pg_1_scedrp) at (8, 2) {};  
        \node[ladstyle] (pg_1_lad) at (7, 3) {};  
        \node[sladstyle] (pg_1_slad) at (3, 7) {};  
        
        \node[scedstyle] (pg_2_sced) at (20, 10) {};  
        \node[scedrpstyle] (pg_2_scedrp) at (20, 12) {};  
        \node[ladstyle] (pg_2_lad) at (20, 13) {};  
        \node[sladstyle] (pg_2_slad) at (20, 15) {};  

        \draw[SCED, line width=0.25pt] (pg_0) -- (pg_1_sced);
        \draw[SCEDRP, line width=0.25pt] (pg_0) -- (pg_1_scedrp);
        \draw[LAD, line width=0.25pt] (pg_0) -- (pg_1_lad);
        \draw[SLAD, line width=0.25pt] (pg_0) -- (pg_1_slad);
        \draw[SCED, line width=0.25pt] (pg_1_sced) -- (pg_2_sced);
        \draw[SCEDRP, line width=0.25pt] (pg_1_scedrp) -- (pg_2_scedrp);
        \draw[LAD, line width=0.25pt] (pg_1_lad) -- (pg_2_lad);
        \draw[SLAD, line width=0.25pt] (pg_1_slad) -- (pg_2_slad);

        \node[scedstyle] at (24, 15) {};
        \draw[SCED, line width=0.25pt] (23, 15) -- (25, 15);
        \node[right] at (25, 15) {\tiny{SCED}};
        
        \node[scedrpstyle] at (24, 12) {};
        \draw[SCEDRP, line width=0.25pt] (23, 12) -- (25, 12);
        \node[right] at (25, 12) {\tiny{SCED+RP}};
        
        \node[ladstyle] at (24, 9) {};
        \draw[LAD, line width=0.25pt] (23, 9) -- (25, 9);
        \node[right] at (25, 9) {\tiny{LAD}};
        
        \node[sladstyle] at (24, 6) {};
        \draw[SLAD, line width=0.25pt] (23, 6) -- (25, 6);
        \node[right] at (25, 6) {\tiny{SLAD}};
        
        \draw[black, line width=0.25pt, dash pattern=on 2pt off 1pt] (23, 3) -- (25, 3);
        \node[right] at (25, 3) {\tiny{Power balance}};
        
        \draw[black, line width=0.25pt, dash pattern=on 2pt off 1pt] (23, 3) -- (25, 3);
        \node[right] at (25, 3) {\tiny{Power balance}};
    \end{tikzpicture}
    }
    
    \caption{%
    Illustration of the active power dispatch obtained by each economic dispatch formulation (\textcolor{SCED}{\textbf{SCED}}/\textcolor{SCEDRP}{\textbf{SCED+RP}}/\textcolor{LAD}{\textbf{LAD}}/\textcolor{SLAD}{\textbf{SLAD}}).
    Both generator's initial output ($t{=}0$, black square) is 0MW.
    At $t \, {=} \, 1$, all formulations meet the 10MW demand.
    At $t \, {=} \, 2$, both SCED and LAD fail to meet the 35MW demand, because generator G2 is not able to ramp sufficiently fast.
    SLAD is able to meet the power balance constraint because it had sufficiently ramped up G2 at $t \, {=} \, 1$.
    }
    \label{fig:example:dispatch:G1G2}
\end{figure}
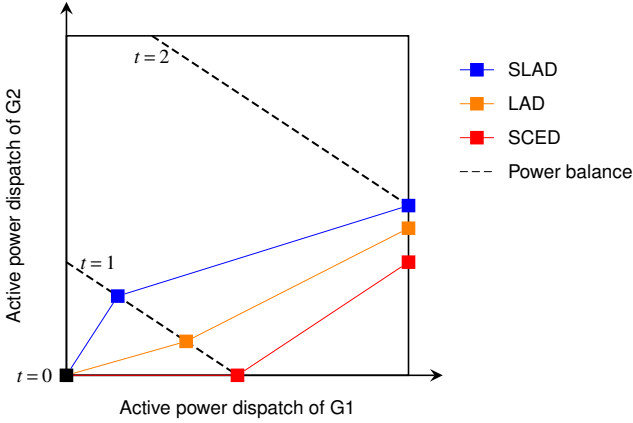

To better understand the behavior of each ED formulation, Figures \ref{fig:example:dispatch} and \ref{fig:example:dispatch:G1G2} illustrate the resulting dispatch decisions at $t \,{=}\, 1, 2$.
Namely, Figure \ref{fig:example:dispatch} illustrates the dispatch of generator G1 and G2 (left and right pane, respectively) at time $t \,{=}\, 0, 1, 2$.
In particular, the figure shows that SLAD and, to a less extent, SCED+RP and LAD, pre-emptively ramp-up generator G2 at $t \, {=} \, 1$, despite the immediate additional cost.
This early ramp-up allows SLAD, SCED+RP and LAD to dispatch generator G2 to a higher value at $t \, {=} \, 2$ than SCED, thereby reducing energy shortage penalty costs.
Nevertheless, unlike SLAD, SCED+RP and LAD still experiences some energy shortages, albeit less than SCED.
This is because the demand forecast used in LAD at $t \, {=} \, 1$ underestimates future demand, thereby resulting in lower dispatch for generator G2 at $t \, {=} \, 1$ than SLAD.
Similarly, the flexible ramp requirement in SCED+RP at $t \, {=} \, 1$ underestimates net load variability, resulting in a lower dispatch for generator G2 at $t \, {=} \, 1$.
Finally, Figure \ref{fig:example:dispatch:G1G2} displays the same solution in the $(p_{1}, p_{2})$ space.
It is evident from Figure \ref{fig:example:dispatch:G1G2} that SCED, SCED+RP and LAD fail to meet the power balance requirement at $t \, {=} \, 2$, and incur energy shortages of 5MW, 3MW and 2MW, respectively.

Naturally, the solution obtained by SCED+RP, and the corresponding costs, heavily depend on the amount of FRP procured by system operators.
In this example, when the FRP requirement at $t \, {=} \, 1$ is no greater than $20$MW, the SCED+RP formulation behaves identically to SCED: this is an indication that the FRP constraint \eqref{eq:stoch:ramp_up_prod} is not binding.
However, increasing the FRP amount from 20MW to 25MW results in a likewise reduction in power shortage at $t \, {=} \, 2$, from 5MW (the same value as SCED) to 0MW.
Further increasing the FRP amount beyond 25MW also eliminates power shortages, but increases operating costs as more expensive generators are dispatched to jointly procure energy and flexible ramping.
Therefore, as noted in \citep{WangHobbs2013_Flexiramp}, identifying the right amount of FRP to be procured is essential to ensure the reliability and economic efficiency of SCED+RP, and an important source of externalities since FRP requirements are set by system operators.

This small numerical example highlights the following lessons.
First, a single-period formulation takes myopic decisions that may fail to account for future ramping needs.
Second, the performance of SCED+RP is sensitive to the amount of FRP being procured, which must be selected a priori to balance reliability and economic costs.
Third, LAD is also sensitive to forecasting errors, especially when future demand is underestimated.
Fourth, SLAD is better at hedging against future uncertainty, without requiring perfect information.

\section{Large-Scale Numerical Experiments}
\label{sec:experiment_RTE}

This section describes the large-scale experiments conducted on the French RTE system with 6708 buses.
Detailed results are reported in Section \ref{sec:results}.

\subsection{Economic Dispatch Formulations}

Table \ref{tab:formulations} summarizes the four economic dispatch formulations considered in the experiment.
As described in Section \ref{sec:formulation}, the SCED, SCED+RP, LAD and SLAD formulations are single-period, single-period with ramping products, multi-period deterministic, and multi-period stochastic formulations, respectively.
In all formulations, time is discretized in 5-minute intervals.
The paper also includes the so-called \emph{perfect dispatch} (PD) \citep{PerfectDispatch}.
The PD formulation takes as input the actual realization of load and renewable production for the entire day, and solves a LAD problem with a whole day 24-hour horizon. This is the best possible operation that can be attained through economic dispatch.
Because it relies on a perfect prediction oracle, it cannot be implemented in practice; instead, PD is used as an \emph{a posteriori} analysis tool \citep{PerfectDispatch}.

\begin{table}[!t]
    \centering
    \caption{The Economic Dispatch Formulations Used in Experiments.}
    \label{tab:formulations}
    \resizebox{\columnwidth}{!}
    {
    \begin{tabular}{ccccccc}
        \toprule
            Problem & Horizon & \# Time steps & \# Scenarios & Source of input data\\
        \midrule
            SCED    & \phantom{0}5 min & 1 & 1 & Current measurements \\
            SCED+RP & \phantom{0}5 min & 1 & 1 & Current measurements \& ramping requirements \\
            LAD     & 60 min & 12 & 1 & Deterministic forecast \\
            SLAD    & 60 min & 12 & 10$^{*}$ & Probabilistic forecast \\
        \midrule
            PD   & 24 hr~~~ & 288 & 1 & Perfect forecast \\
        \bottomrule
    \end{tabular}
    }\\
    \footnotesize{All formulations use 5-minute time steps. 
    \\ $^{*}$SLAD considers uncertainty captured by 10 scenarios each for load, wind, and solar.}
\end{table}

\subsection{Experimental Setting}
\label{sec:experiment_RTE:setting}

\begin{figure}[!t]
    \centering
    \includegraphics[width=0.9\columnwidth]{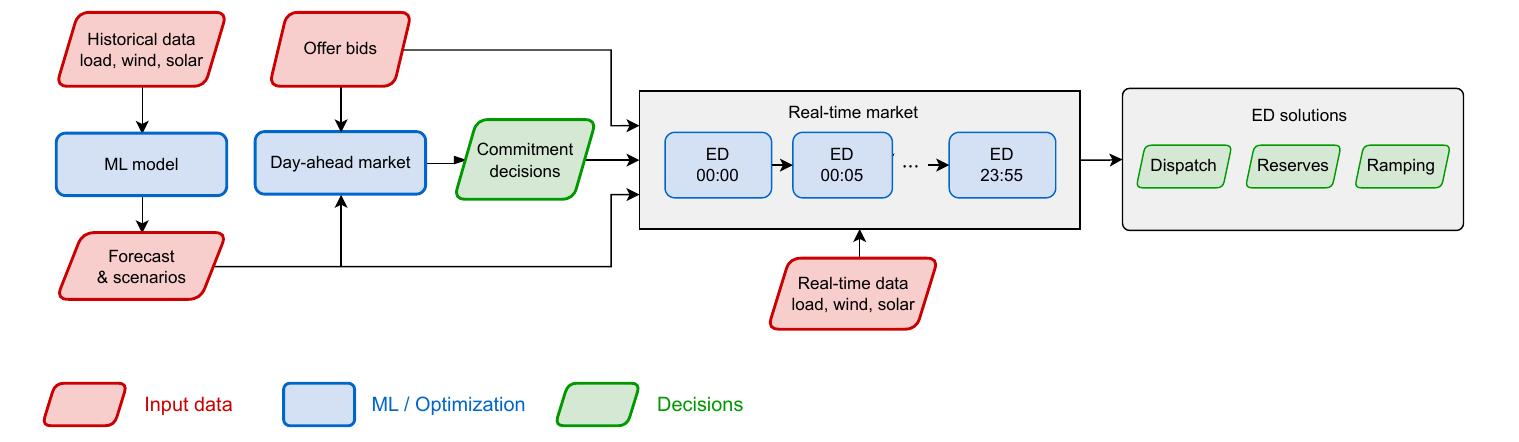}
    \caption{The Simulation Pipeline.}
    \label{fig:simulation_pipeline}
\end{figure}

The performance of SLAD is evaluated on the over 6,000 buses RTE system for the entire year in 2018.
Namely, the paper executes a simulation for every day of the year, using a realistic pipeline depicted in Figure \ref{fig:simulation_pipeline}, which reproduces the interaction between day-ahead and real-time markets in the US.
The detailed setup of the case study is as follows.
The real-time market pipeline is executed over the day, i.e., an economic dispatch problem is solved every 5 minutes, and it outputs the energy dispatch and reserves for each generator which become inputs for the following period.

The time series for load and renewable generation, as well as economic data for the generators, were reconstructed following the disaggregation methods in \citep{chatzos2022data}.
The case study sets high penalty costs with load shedding at \$100,000/MWh, operating reserve shortages at \$50,000/MWh, regulating reserve shortage at \$55,000/MWh, transmission violation cost \$1,500/MWh, ramping capability shortage at \$30/MWh (which is the price used by some ISOs). 
Artificial generators that represent imports have a high capacity and a cost of \$1000/MWh.
Therefore, the system should not encounter load shedding or reserve shortages as import generators provide sufficient high-cost energy if needed.

All optimization problems are formulated in JuMP \citep{DunningHuchetteLubin2017_JuMP} and solved using Gurobi 10.0 with 10 threads.
The SCED, SCED+RP and LAD problems employ Gurobi's barrier method (method parameter set to 2).
SLAD’s accelerated Benders algorithm is implemented in Julia, solving the BD subproblems in parallel—each with one core using the barrier method—and the master problem with 10 threads.
Consequently, SLAD uses the same total of 10 threads as SCED, SCED+RP, and LAD.
The BD algorithm uses a maximum number of iterations of 100, and a relative gap tolerance of $\epsilon \, {=} \, 10^{-5}$. 
Experiments are carried out on the Phoenix cluster at the Georgia Institute of Technology, Atlanta, Georgia \citep{PACE}.

\subsection{Forecast Generation}
\label{sec:experiment_RTE:forecast}

Deep learning models were used to generate scenarios and forecasts for the RTE system \citep{scenariogeneration}.
The DeepAR model \citep{salinas2020deepar} is used to generate wind scenarios, and a Temporal Fusion Transformer model
\citep{lim2021temporal} is used for solar and load scenario generation.
The DeepAR and TFT models are probabilistic forecasting models, from whose prediction intervals one can sample multiple scenarios. To build confidence-aware prediction sets, we may also leverage more advanced techniques such as generative models \citep{colombo2024normalizing} or conformal prediction \citep{xu2024conformal,xu2023sequential} in the future.
For LAD, a deterministic forecast is obtained by taking the average of 500 random scenarios.
For SLAD, 10 scenarios are considered.
The interplay between the quality of the forecasts/scenarios and the performance of LAD and SLAD is further analyzed in Section \ref{sec:sensitivity:scenarios}.

\subsection{Cost Comparison}
\label{sec:experiment_RTE:costs}

The paper uses the following methodology to evaluate costs and compare the performance of each ED formulation.
First, at each time step of a simulation, operating costs are evaluated from generators' energy and reserve dispatch, and import costs.
Then, constraint violation penalty costs are computed, which include reserve shortage costs, ramping products shortage costs, transmission violation penalties and energy shortage/surplus penalties.
This yields a combination of operating and penalty costs at each time step, which are then aggregated over the day.
Finally, recall that each simulation is executed once per ED formulation, which is achieved by using the same random seed when initializing the simulation environment.
This ensures that each ED formulation is executed against identical realizations of load and renewable output, thereby enabling a fair comparison.

\section{Result Analysis}
\label{sec:results}

This section presents the results of the large-scale experiment on the RTE system.
Section \ref{sec:results:computimg_time} reports the computational performance of each formulation.
Section \ref{sec:results:overview} compares the economic performance of the various ED formulations, which is further detailed in Sections \ref{sec:results:multi-period} and \ref{sec:results:stochastic}.
Finally, Section \ref{sec:results:detailed} provides a detailed analysis of the behavior of SLAD.

\subsection{Computing Times}
\label{sec:results:computimg_time}

\begin{figure}[!t]
    \centering
    \includegraphics[width=0.6\columnwidth]{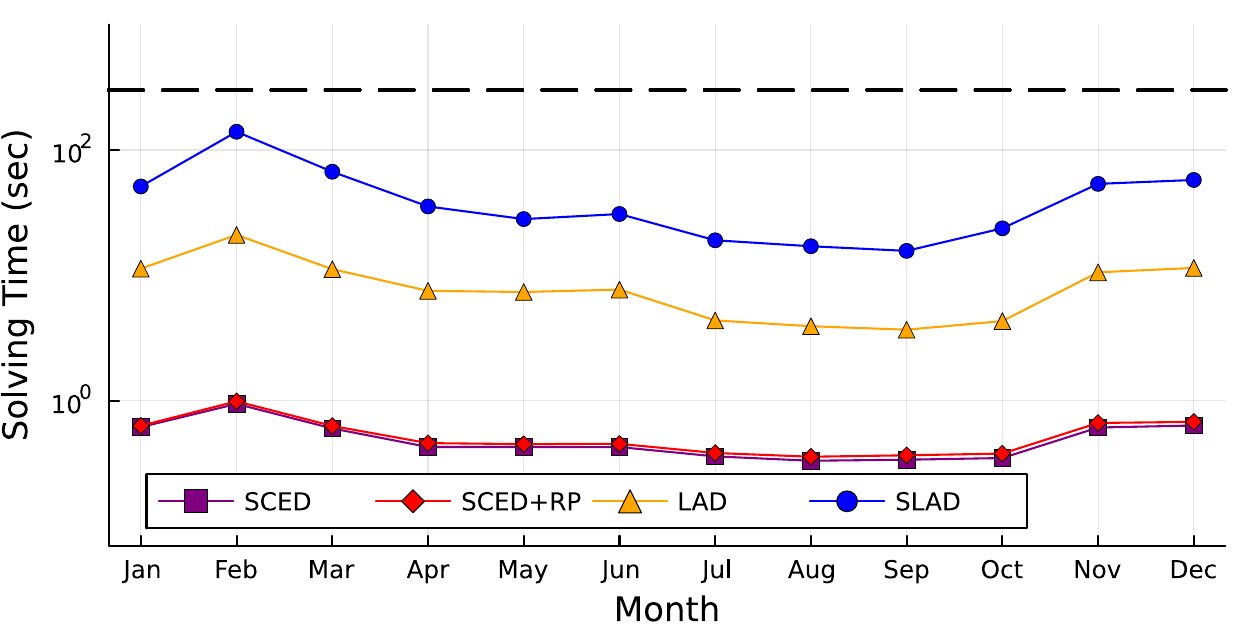}
    \caption{
    Average solution times of ED formulations of each month.
    The dashed \textbf{black} line denotes the 5-minute frequency of real-time market clearing.}
    \label{fig:sol_time}
\end{figure}

Figure \ref{fig:sol_time} reports the solution time of each formulation over the year.
For each formulation, the figure displays the average computing time for each month.
First, observe that computing times are higher during winter, and lower during summer.
Indeed, demand is highest during the winter, thereby increasing network congestion and the number of PTDF constraints that must be added to the ED formulations.
This naturally increases memory requirements and computing times.
Second, most instances are solved within the 5-minute time frame of real-time operations.
SCED and SCED+RP instances are typically solved within one second, and LAD instances are solved within 3 to 22 seconds, on average.
Similarly, although they require longer computing times, SLAD instances are solved within an average of 15 to 140 seconds.
In fact, 99.77\% of SLAD instances reach the prescribed $0.001\%$ optimality gap tolerance within 5 minutes.
These results confirm that multi-period ED formulations, both deterministic and stochastic, are computationally tractable for real-time market operations.

\subsection{Results Overview}
\label{sec:results:overview}

Table \ref{tab:total_cost_table:DL} reports, for each ED
formulation, the average daily costs and savings over the year.
Results are reported for low-load days and high-load days, which are
defined as days where the average load is less than 50GW and greater
than 50GW, respectively.  
Note that low- and high-load days follow a seasonal trend, with overall demand being higher in winter and lower in summer.
The last three columns in Table \ref{tab:total_cost_table:DL} report results over the whole year.
For day $d$, the daily savings of formulation X compared to SCED are
given by
    \begin{align}
        \frac{C^{\text{SCED}}_{d} - C^{\text{X}}_{d}}{C^{\text{SCED}}_{d}},
    \end{align}
where $C^{\text{X}}_{d}$ (resp. $C^{\text{SCED}}_{d}$) denotes the
total costs (including penalties) attained by formulation X
(resp. SCED) on day $d$; positive values indicate net savings.  For
instance, SCED+RP achieves an average daily
cost of 13.09M\$ on low-load days, 34.56M\$ on high-load days, and
23.62M\$ over the year.  The corresponding savings relative to SCED
are 0.48\% on low-load days, 0.78\% on high-load days, and 0.70\% over
the whole year.
Specifically, SLAD achieves 65.7\% more savings overall compared to SCED+RP and further enhances performance over LAD by incorporating multiple scenarios, rather than depending on a single deterministic forecast.

\begin{table*}[!t]
    \centering
    \caption{
        Comparison of daily operating costs in RTE over 2018.
        Low (resp. high) load days are days where the average load is less (resp. greater than) 50GW.
        All costs are in M\$ per day.
        Savings are relative to SCED.
    }
    \label{tab:total_cost_table:DL}
    \resizebox{\columnwidth}{!}
    {
    \begin{tabular}{llrrrrrrrr}
        \toprule
            & \multicolumn{3}{c}{Low load} 
            & \multicolumn{3}{c}{High load}
            & \multicolumn{3}{c}{All}\\
        \cmidrule(llr){2-4}
        \cmidrule(llr){5-7}
        \cmidrule(llr){8-10}
        Formulation 
            & Cost
            & Savings
            & Savings(\%)
            & Cost
            & Savings
            & Savings(\%)
            & Cost
            & Savings
            & Savings(\%)\\
        \midrule
        SCED
           & 13.15 & 0.00 & -- 
           & 34.84 & 0.00 & -- 
           & 23.79 & 0.00 & -- 
           \\
        SCED+RP
           & 13.09 & 0.06 & 0.48 \% 
           & 34.56 & 0.27 & 0.78 \% 
           & 23.62 & 0.17 & 0.70 \% 
           \\
        LAD
           & 13.02 & 0.14 & 1.04 \% 
           & 34.47 & 0.36 & 1.04 \% 
           & 23.54 & 0.25 & 1.04 \% 
           \\
        SLAD
           & \textbf{13.01} & \textbf{0.15} & \textbf{1.11} \% 
           & \textbf{34.43} & \textbf{0.41} & \textbf{1.17} \% 
           & \textbf{23.51} & \textbf{0.28} & \textbf{1.16} \% 
           \\
        \midrule
        PD$^{\dagger}$
           & 12.98 & 0.17 & 1.32 \% 
           & 34.26 & 0.57 & 1.64 \% 
           & 23.42 & 0.37 & 1.55 \% 
           \\
        \bottomrule
    \end{tabular}
    }\\
    \footnotesize{$^{\dagger}$PD is not practical for real-world application as it presupposes an accurate forecast over the next 24 hours. This is the best possible operation that can be attained through economic dispatch.}
\end{table*}

\begin{figure}[!t]
    \centering
    \includegraphics[width=0.48\columnwidth]{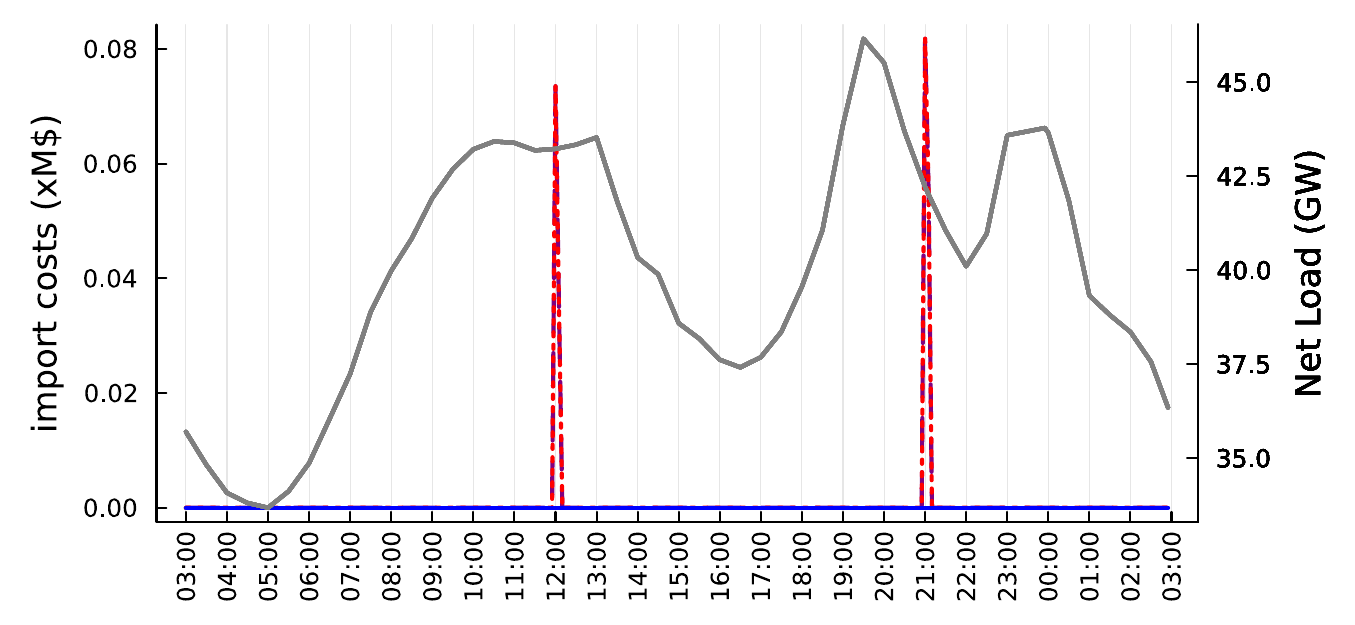}
    \hfill
    \includegraphics[width=0.48\columnwidth]{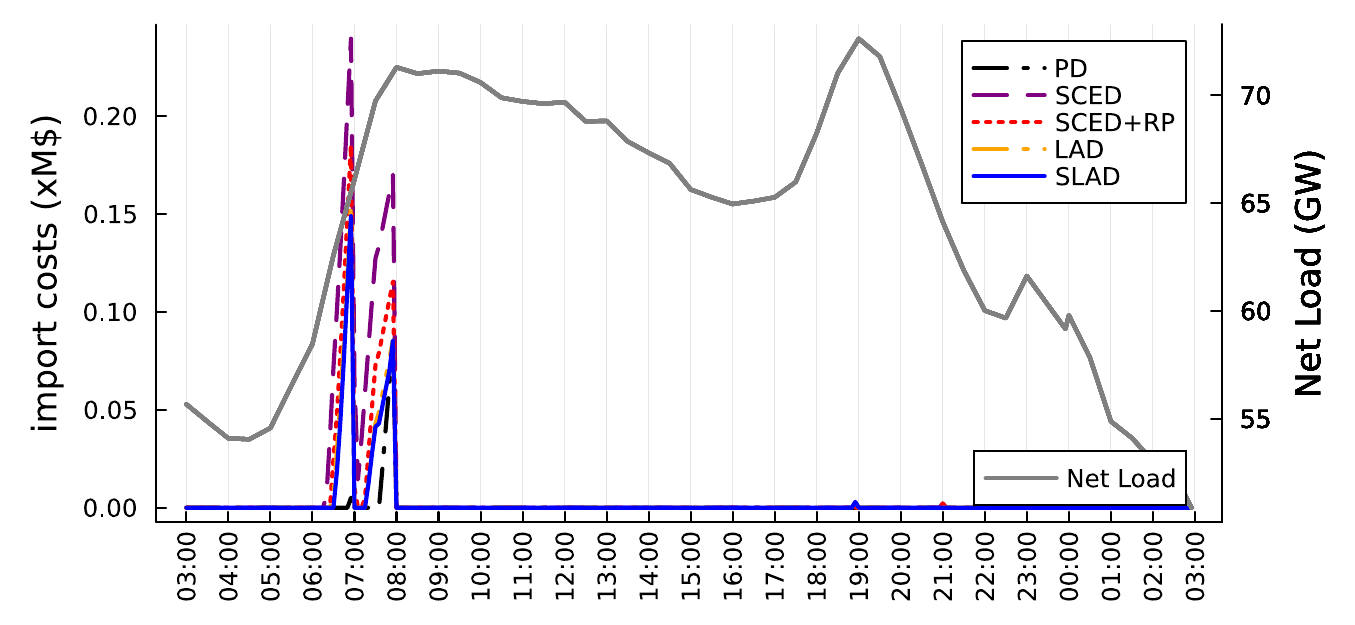}
    \caption{Import cost comparison on medium-load (10/20) and high-load (01/30) days. The import cost on the low-load day (07/25) is zero, so the corresponding graph is omitted here.}
    \label{fig:results:itemized_costs}
\end{figure}

Table \ref{tab:cost_comparison:three_cases} provides an itemized cost
comparison between all ED formulations, on three selected days: July 25th, October 20th, and January 30th.  For each day and formulation,
Table \ref{tab:cost_comparison:three_cases} reports the total cost of
domestic energy production (Energy), total import costs (Import), and
the sum of no-load costs, reserve procurement costs and constraint
violation penalties (Other).  Relative savings compared to SCED are
also reported.  Figure \ref{fig:results:itemized_costs} depicts the evolution of import costs with the right-axis depicts the net load in GW over time.
    
July 25th is a summer day with low demand; the system has ample
generation capacity and flexibility, and all formulations achieve
almost the same cost.  October 20th is a Fall day during which the
system experiences moderate net load ramping, especially in the
morning and evening hours, thereby requiring more flexibility from the
conventional generation fleet.  January 30th is a Winter day with high
net load and high ramping needs in the morning.  This day also
experiences higher-than-usual errors in the day-ahead forecast for
load and renewables, which results in insufficient commitments and
capacity shortages throughout the day.  A detailed analysis of this
day is presented in Section \ref{sec:results:detailed}.

\begin{table}[!t]
    \centering
    \caption{Itemized daily cost comparison for three representative days in 2018. Cost items and absolute savings relative to SCED are presented in millions of dollars (M\$), with savings also expressed as a percentage (\%). }
    \label{tab:cost_comparison:three_cases}
    \resizebox{0.7\columnwidth}{!}{
    \begin{tabular}{llrrrr|r}
        \toprule
        Date
        & Cost item
           & \multicolumn{1}{c}{SCED\phantom{\%}}
           & \multicolumn{1}{c}{SCED+RP\phantom{\%}}
           & \multicolumn{1}{c}{LAD\phantom{\%}}
           & \multicolumn{1}{c}{SLAD\phantom{\%}}
           & \multicolumn{1}{c}{PD$^{\dagger}$\phantom{\%}}
        \\
        \midrule
        07/25
           & Energy
              & 9.74 \phantom{\%} 
              & 9.74 \phantom{\%} 
              & 9.74 \phantom{\%} 
              & 9.74 \phantom{\%} 
              & 9.74 \phantom{\%} 
              \\
           & Import
              & 0.00 \phantom{\%} 
              & 0.00 \phantom{\%} 
              & 0.00 \phantom{\%} 
              & 0.00 \phantom{\%} 
              & 0.00 \phantom{\%} 
              \\
           & Other$^{\dagger}$
              & 5.13 \phantom{\%} 
              & 5.13 \phantom{\%} 
              & 5.13 \phantom{\%} 
              & 5.13 \phantom{\%} 
              & 5.13 \phantom{\%} 
              \\
        \cmidrule{2-7}
           & Total
              & 14.87 \phantom{\%} 
              & 14.87 \phantom{\%} 
              & 14.87 \phantom{\%} 
              & 14.87 \phantom{\%} 
              & 14.87 \phantom{\%} 
              \\
           & Savings
              & 0.00 \phantom{\%} 
              & 0.00 \phantom{\%} 
              & 0.00 \phantom{\%} 
              & 0.00 \phantom{\%} 
              & 0.00 \phantom{\%} 
              \\
           & Savings(\%)
              & 0.00 \%
              & -0.01 \% 
              & \textbf{0.00} \% 
              & \textbf{0.00} \% 
              & 0.01 \% 
              \\
        \midrule
        10/20
           & Energy
              & 8.64 \phantom{\%} 
              & 8.66 \phantom{\%} 
              & 8.65 \phantom{\%} 
              & 8.66 \phantom{\%} 
              & 8.65 \phantom{\%} 
              \\
           & Import
              & 0.25 \phantom{\%} 
              & 0.25 \phantom{\%} 
              & 0.00 \phantom{\%} 
              & 0.00 \phantom{\%} 
              & 0.00 \phantom{\%} 
              \\
           & Other$^{\dagger}$
              & 5.11 \phantom{\%} 
              & 5.10 \phantom{\%} 
              & 5.12 \phantom{\%} 
              & 5.11 \phantom{\%} 
              & 5.11 \phantom{\%} 
              \\
        \cmidrule{2-7}
           & Total
              & 14.00 \phantom{\%} 
              & 14.01 \phantom{\%} 
              & 13.77 \phantom{\%} 
              & 13.77 \phantom{\%} 
              & 13.76 \phantom{\%} 
              \\
           & Savings
              & 0.00 \phantom{\%} 
              & -0.01 \phantom{\%} 
              & 0.23 \phantom{\%} 
              & 0.23 \phantom{\%} 
              & 0.24 \phantom{\%} 
              \\
           & Savings(\%)
              & 0.00 \%
              & -0.08 \% 
              & \textbf{1.70} \% 
              & 1.69 \% 
              & 1.71 \% 
              \\
        \midrule
        01/30
           & Energy
              & 16.71 \phantom{\%} 
              & 16.73 \phantom{\%} 
              & 16.74 \phantom{\%} 
              & 16.75 \phantom{\%} 
              & 16.76 \phantom{\%} 
              \\
           & Import
              & 2.23 \phantom{\%} 
              & 1.27 \phantom{\%} 
              & 0.87 \phantom{\%} 
              & 0.80 \phantom{\%} 
              & 0.22 \phantom{\%} 
              \\
           & Other$^{\dagger}$
              & 10.32 \phantom{\%} 
              & 10.31 \phantom{\%} 
              & 10.29 \phantom{\%} 
              & 10.28 \phantom{\%} 
              & 10.26 \phantom{\%} 
              \\
        \cmidrule{2-7}
           & Total
              & 29.26 \phantom{\%} 
              & 28.31 \phantom{\%} 
              & 27.90 \phantom{\%} 
              & 27.83 \phantom{\%} 
              & 27.24 \phantom{\%} 
              \\
           & Savings
              & 0.00 \phantom{\%} 
              & 0.95 \phantom{\%} 
              & 1.36 \phantom{\%} 
              & 1.43 \phantom{\%} 
              & 2.02 \phantom{\%} 
              \\
           & Savings(\%)
              & 0.00 \%
              & 3.25 \% 
              & 4.64 \% 
              & \textbf{4.89} \% 
              & 6.92 \% 
              \\
        \bottomrule
    \end{tabular}
    }
    \\
    \footnotesize{
    $^{\dagger}$PD is not practical for real-world application as it presupposes an accurate forecast over the next 24 hours. This is the best possible operation that can be attained through economic dispatch.
    }
\end{table}

\subsection{Individual Factor Analysis}
The next sections further analyze the ED formulations by focusing on individual factors. Section \ref{sec:results:ramping-products} examines the effect of incorporating ramping products, Section \ref{sec:results:multi-period} compares single-period versus multi-period formulations, and Section \ref{sec:results:stochastic} explores the differences between deterministic and stochastic approaches. These analyses provide a deeper understanding of the elements that influence performance.

\subsubsection{Without vs With Ramping Products}
\label{sec:results:ramping-products}

This section examines the impact of incorporating ramping products by comparing the SCED model with the SCED that includes ramping products (SCED+RP). As shown in the first two rows of Table \ref{tab:total_cost_table:DL}, SCED+RP achieves a cost savings of 0.48\% under low load conditions, increasing to 0.78\% under high load conditions. On average, incorporating ramping products leads to a 0.70\% reduction in costs. These results indicate that, by incorporating ramping products, the system gains additional flexibility to manage net load variability more effectively. The savings are consistent and most significant under high-load conditions, where the need for flexibility is greatest.

\subsubsection{Single-Period vs Multi-Period}
\label{sec:results:multi-period}

It is interesting to understand the benefits of using a deterministic multi-period economic dispatch formulation, such as LAD, instead of a deterministic single-period formulation.
Note that, as illustrated in the example of Section \ref{sec:example}, the performance of LAD is affected by the accuracy of the forecast it receives as input; this aspect is further analyzed in Section \ref{sec:sensitivity:scenarios}. Over the year, as noted in Table \ref{tab:total_cost_table:DL}, LAD saves an average of 1.04\% in daily costs compared to SCED.

Table \ref{tab:cost_comparison:three_cases} and Figure \ref{fig:results:itemized_costs} enable a more granular comparison.
On July 25th, there is no significant difference in cost between single- and multi-period formulations.
On October 20th, the system experiences higher net load ramping, especially in the morning and evening hours, thereby requiring more flexibility from the generation fleet.
While  LAD is able to accommodate the ramping needs, and achieves the same cost as PD, SCED needs to import energy twice during the day, which results in higher total costs.
Finally, recall that, on January 30th, because of errors in the day-ahead forecast used in the unit commitment, the online generation capacity is not enough to meet all energy and reserve requirements.
As a result, all formulations import energy during the morning, but again SCED performs the worst, with LAD saving on imports and 4.64\% overall compared to SCED.
These results demonstrate that the multi-period formulations, compared to the ramping products, better anticipate future ramping needs which, although more expensive in the very-short term, yields savings in the long term.

\subsubsection{Deterministic vs Stochastic}
\label{sec:results:stochastic}

Consider now the potential benefits of a stochastic multi-period
economic dispatch (SLAD) formulation instead of a deterministic
multi-period formulation (LAD).  Table
\ref{tab:total_cost_table:DL} shows that SLAD saves an
average of 1.16\% per day in total costs.  Savings are
seasonal, with SLAD saving an average 1.17\% on high-load
days, and an average 1.11\% on low-load days.  These observations
are further supported by the results in Table
\ref{tab:cost_comparison:three_cases}.  On July 25th and October
20th, SLAD and LAD achieve about the same costs, and match the
performance of PD.  However, on January 30th, which experiences
higher ramping needs, SLAD reduces import costs by 5\% compared to
LAD, and total costs by 0.25\%.  Over the year, the savings
achieved by SLAD are higher than those of LAD by 12\%
overall, 14\% on high-load days, and 7\% on low-load
days.  This difference demonstrates the benefits of considering
multiple scenarios instead of a single point forecast: the former
hedges against uncertainty whereas the latter is more sensitive to
forecast errors.

\subsection{Detailed Analysis}
\label{sec:results:detailed}

This section further investigates the behavior of SCED+RP and SLAD on Jan
30th, a Winter day with a high net load.  As shown in Table
\ref{tab:cost_comparison:three_cases} and Figure
\ref{fig:results:itemized_costs}, the main difference in costs stems
from imports, with SCED+RP incurring 60\% more imports between 6am and
8am.  Recall that, because all simulations use the same commitment
decisions and load/renewable profiles, the differences in costs are
purely a consequence of how dispatch decisions are made throughout the
day.  It is therefore important to understand the underlying factors
explaining the difference in behavior between SCED+RP and SLAD.  The rest
of this section focuses on the morning hours between 5am and 8am,
which is when most differences occur.

\begin{figure}[!t]
    \centering
    \includegraphics[width=0.48\columnwidth]{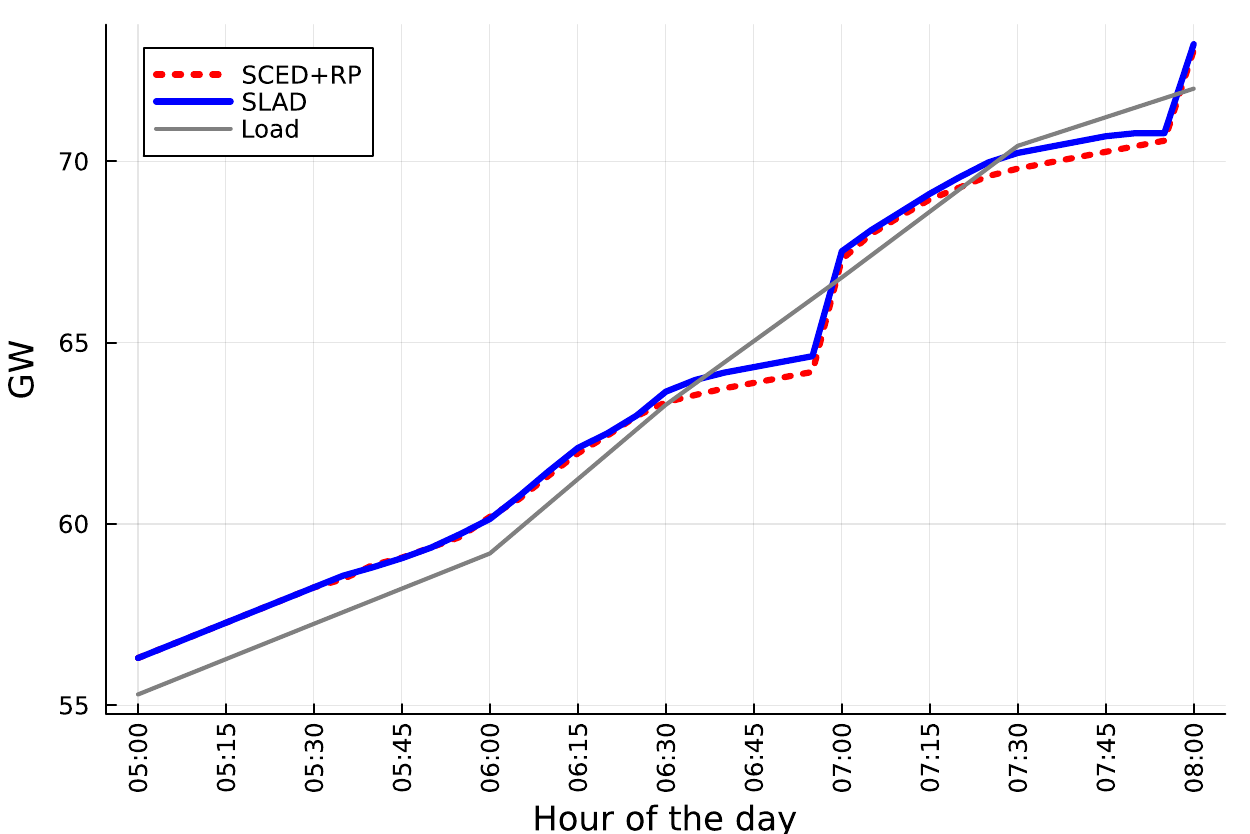}
    \hfill
    \includegraphics[width=0.48\columnwidth]{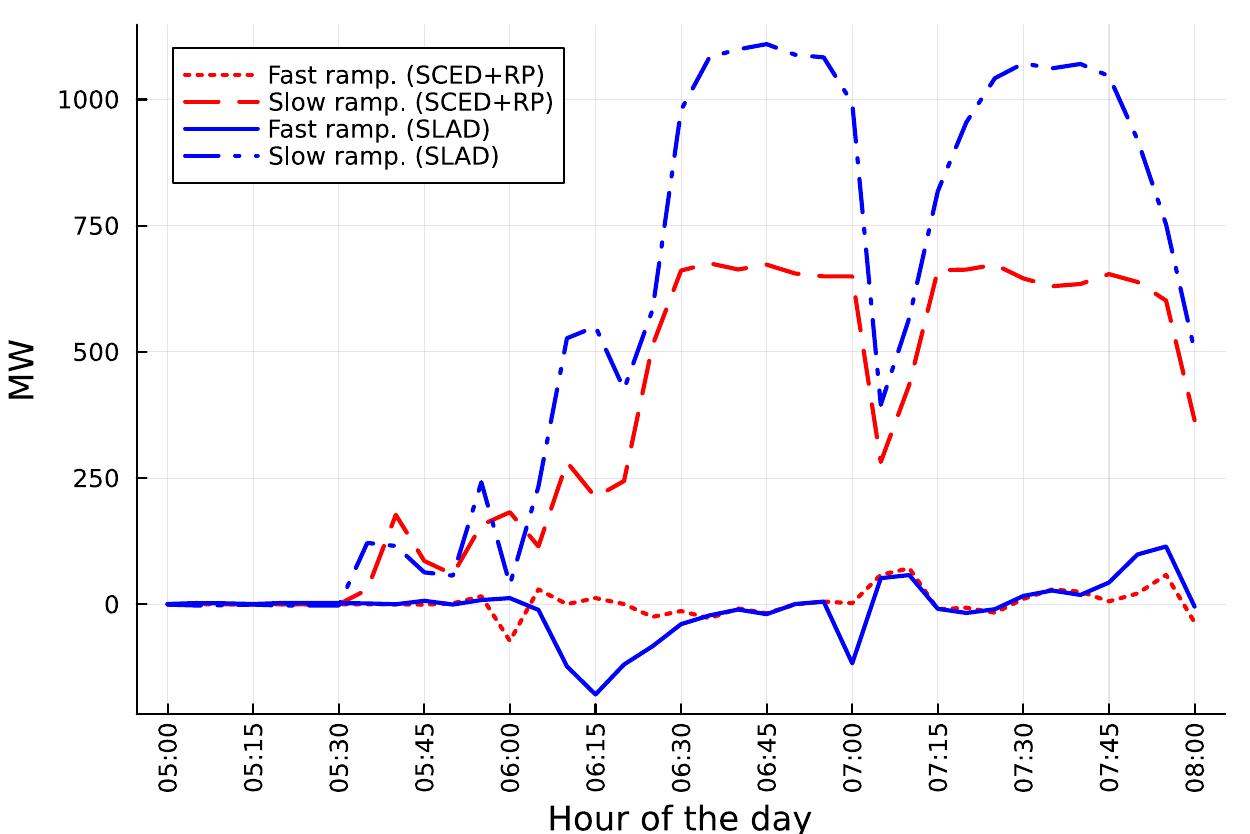}
    \caption{
        Left: Total available capacity of SCED+RP and SLAD between 5am and 8am on January 30th, 2018.\\
        Right: Difference in available capacity between SCED+RP and SCED, SLAD and SCED for slow-ramping and fast-ramping generators, on January 30th, 2018. Slow-ramping (resp. fast-ramping) generators have a ramping rate below (resp. above) 1\% of their total capacity per minute. Positive (resp. negative) numbers indicate the available capacity in SCED+RP or SLAD is higher (Resp. lower) than in SCED.}
    \label{fig:2018-01-30:analysis:avai_cap}
\end{figure}

The available capacity of generator $g$ at time $t$ is denoted by $\pg^{\text{avail}}_{g, t}$.
It is defined as the maximum amount of power this generator can produce, given its operating parameters and previous (5 minutes earlier) dispatch level.
Namely,
\begin{align}
    \label{eq:results:available_capacity}
    \pg^{\text{avail}}_{g, t} &= \min \left\{ \pgmax_{g, t}, \pg_{g,t-1} + \rrup_{g,t} \dt \right\},
\end{align}
where $\pgmax_{g, t}, \pg_{g,t-1}$ and $\rrup_{g,t}$ are the generator's maximum output at time $t$, previous dispatch setpoint, and ramping rate, respectively.
Note that a generator's dispatch cannot exceed its available capacity.
Therefore, if the load exceeds the total available capacity, there must be imports to balance the system.
Also note that imports may still occur even when the total available capacity exceeds loads, in order to satisfy reserve requirements and/or transmission constraints.

The left figure in Figure \ref{fig:2018-01-30:analysis:avai_cap} displays the total
available capacity of SCED+RP and SLAD between 5am and 8am, as
well as the total demand.  There are two important observations from
the figure.  On the one hand, the total load exceeds the available
capacity of SCED+RP from 6:30am to 8am, except for a short period of
time just after 7am. This corresponds exactly to the window during
which SCED+RP imports, thus validating that available capacity is a
useful indicator of the system reliability.  On the one hand, SLAD
always has a higher available capacity than SCED+RP.  The difference
is first noticeable at 5:30am, about one hour before SCED starts
importing.  This delay matches the planning horizon of SLAD, which
suggests that SLAD is able to anticipate the possibility of
(costly) imports at 6:20am, and therefore takes preemptive action.
Nevertheless, SLAD still incurs imports between 6:30am and 7am,
albeit a lot less than SCED+RP, because of insufficient
online generation capacity at that time.  This demonstrates the
importance of intra-day commitments for dealing with unforeseen conditions.

To better understand the above results, the right figure in Figure
\ref{fig:2018-01-30:analysis:avai_cap} displays the difference in available capacity between SLAD and SCED, SCED+RP and SCED over time,
disaggregated between slow- and fast-ramping generators.
Slow-ramping (resp. fast-ramping) generators are defined as generators whose
ramping rate is less than (resp. greater than) 1\% of their total
capacity per minute.
It is clear that
the difference between SLAD and SCED+RP mostly stems from slow-ramping
generators.  The difference in behavior is most apparent after 5:30am,
echoing the results of the left figure in Figure \ref{fig:2018-01-30:analysis:avai_cap}.
Finally, the sharp drop in the difference between SLAD and SCED+RP around
7am and 8am is explained by generators starting up at those times.

The sharp difference in available capacity of slow-ramping generators
seen in the right figure in Figure \ref{fig:2018-01-30:analysis:avai_cap}
can be explained by Equation \eqref{eq:results:available_capacity}.
Indeed, in the RTE system, slow-ramping generators are mostly nuclear
generators with a large (about 1GW) maximum capacity but slow ramping
rates.  Thus, their available capacity is constrained by ramping,
i.e., $\pg^{\text{avail}}_{g, t} \, {=} \, \pg_{g,t-1} \, {+} \,
\rrup_{g,t} \dt$ in Eq. \eqref{eq:results:available_capacity}.  In
contrast, for fast-ramping generators, the relation
$\pg^{\text{avail}}_{g, t} \, {=} \, \pgmax_{g, t}$
typically holds.  Therefore, at time $t$, the total available capacity
is roughly equal to
    \begin{align}
        \label{eq:results:available_capacity:total}
        \pg^{\text{avail}}_{t} &= \sum_{g \in \mathcal{G}_{s}} \left( \pg_{g,t-1} +  \rrup_{g,t} \dt \right) + \sum_{g \in \mathcal{G}_{f}} \pgmax_{g, t},
    \end{align}
where $\mathcal{G}_{f}$ and $\mathcal{G}_{s}$ are slow and fast-ramping generators, respectively.

\begin{figure}[!t]
    \centering
    \includegraphics[width=0.48\columnwidth]{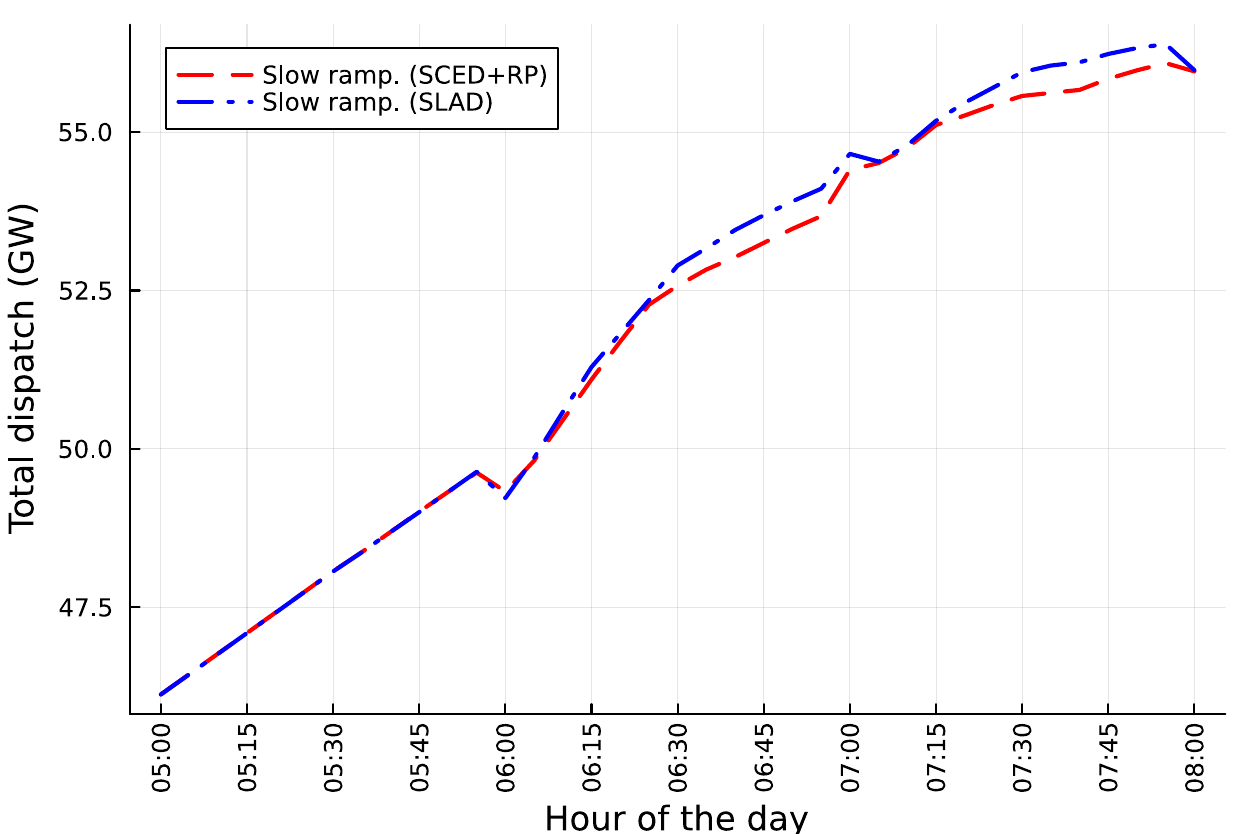}
    \hfill
    \includegraphics[width=0.48\columnwidth]{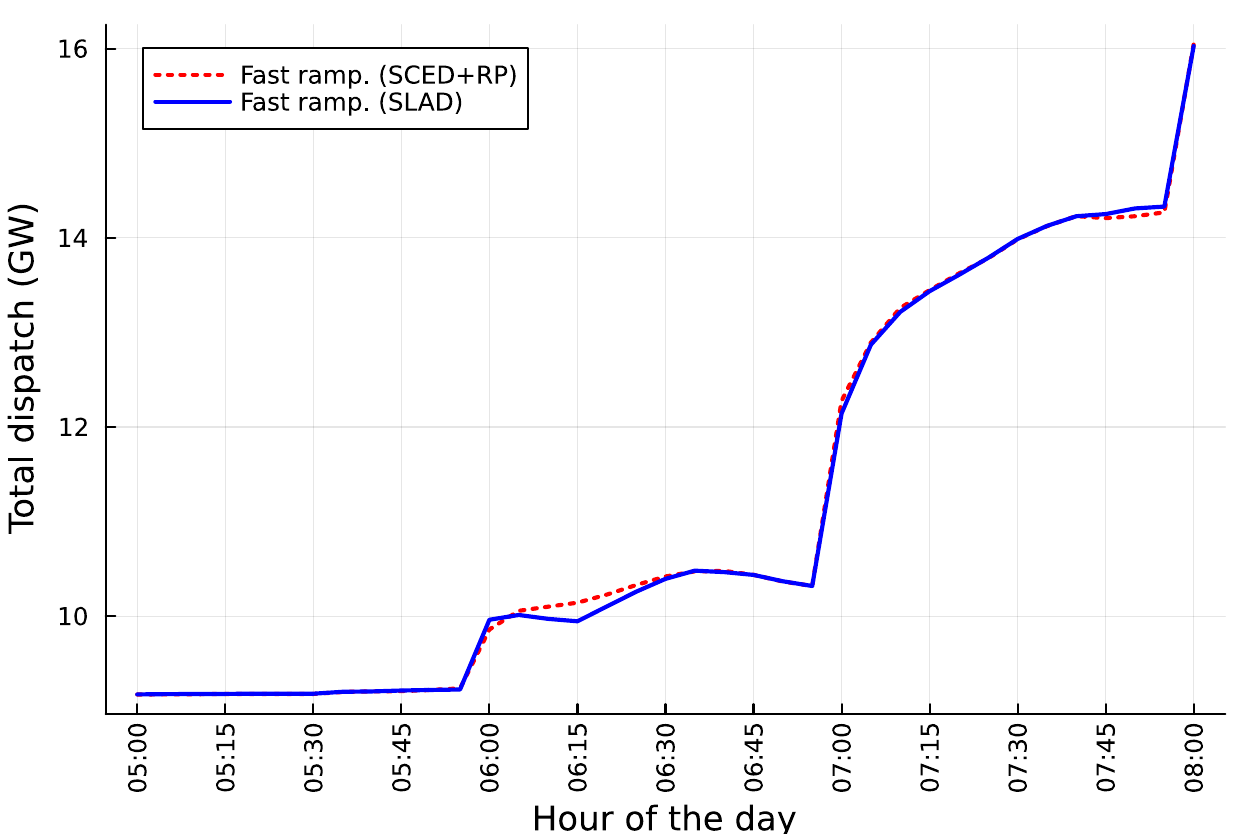}
    \caption{Total dispatch comparison between between SCED+RP and SLAD for slow-ramping (left) and fast-ramping (right) generators between 5:00 and 8:00 on January 30th, 2018}
    \label{fig:2018-01-30:analysis:total_dispatch}
\end{figure}

From Equation \eqref{eq:results:available_capacity:total}, it becomes
evident that, to increase the total available capacity at time $t$,
one must either (i) commit additional units, or (ii) preemptively
increase the dispatch of slow-ramping generators.  The former is
achieved by intra-day commitments, which is out of the scope of the ED
formulations considered in this work.  For the latter, Figure
\ref{fig:2018-01-30:analysis:total_dispatch} displays the total dispatch of
slow- and fast-ramping generators in SLAD and SCED+RP throughout
the study period.
For slow-ramping generators, there is a clear
difference between SLAD and SCED+RP, with the total dispatch of
slow-ramping generators higher in SLAD than in SCED+RP.  

\begin{figure}[!t]
    \centering
    \includegraphics[width=0.5\columnwidth]{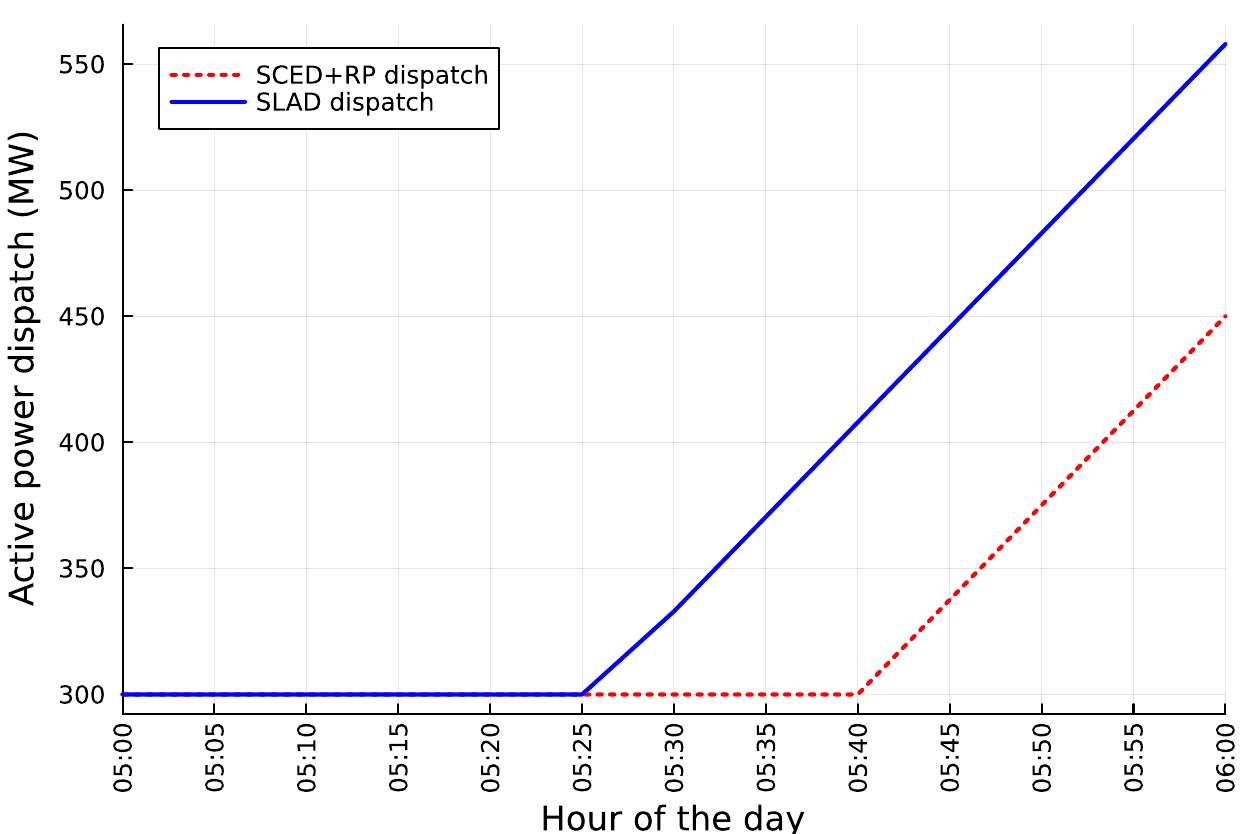}
    \caption{Dispatch of nuclear generator \texttt{CHOO2N01} from SCED+RP and SLAD solutions}
    \label{fig:2018-01-30:analysis:dispatch:CHOO2N01}
\end{figure}

Figure
\ref{fig:2018-01-30:analysis:dispatch:CHOO2N01} illustrates SLAD
and SCED+RP dispatches of nuclear generator \texttt{CHOO2N01} between 5am and 6am.
The figure shows that SLAD starts ramping up \texttt{CHOO2N01} as early as 5:25am, roughly one hour before imports incur, which corresponds to SLAD planning horizon. 
SLAD foresees the possibility of having to rely on imports around 6:30am.
To avoid these costly imports, SLAD preemptively ramps up several nuclear generators, thereby ensuring enough domestic production in the subsequent hours. This preemptive ramping up of nuclear generators enables SLAD to reduce the duration and volume of imports compared to SCED+RP, thereby realizing substantial savings.
SCED+RP starts ramping up the slow-ramping rate generators later than SLAD, partly because the ramping products are designed for only 30 minutes and 10 minutes, durations that are insufficient to adequately accommodate future ramping needs.
Several other nuclear generators display a similar behavior, which accounts for most of the difference in total dispatch of slow-ramping generators between SLAD and SCED+RP.

It is worth noting that SCED+RP, as a single-period formulation, cannot take future commitment decisions into account.
Therefore, for example, even if a cheap generator is online in the future, it will not be able to be used to procure current ramping requirements in SCED+RP.
However, this information can be foreseen by the multi-period formulation SLAD, which is able to determine whether the ramping needs will be satisfied with the generators coming online.

\section{Sensitivity Analysis}
\label{sec:sensitivity}

This section further analyzes the sensitivity of SLAD to the input probabilistic forecast, and to the length of the look-ahead horizon.

\subsection{Impact of the probabilistic forecast}
\label{sec:sensitivity:scenarios}

\begin{figure}[!t]
    \centering
    \includegraphics[width=0.48\columnwidth]{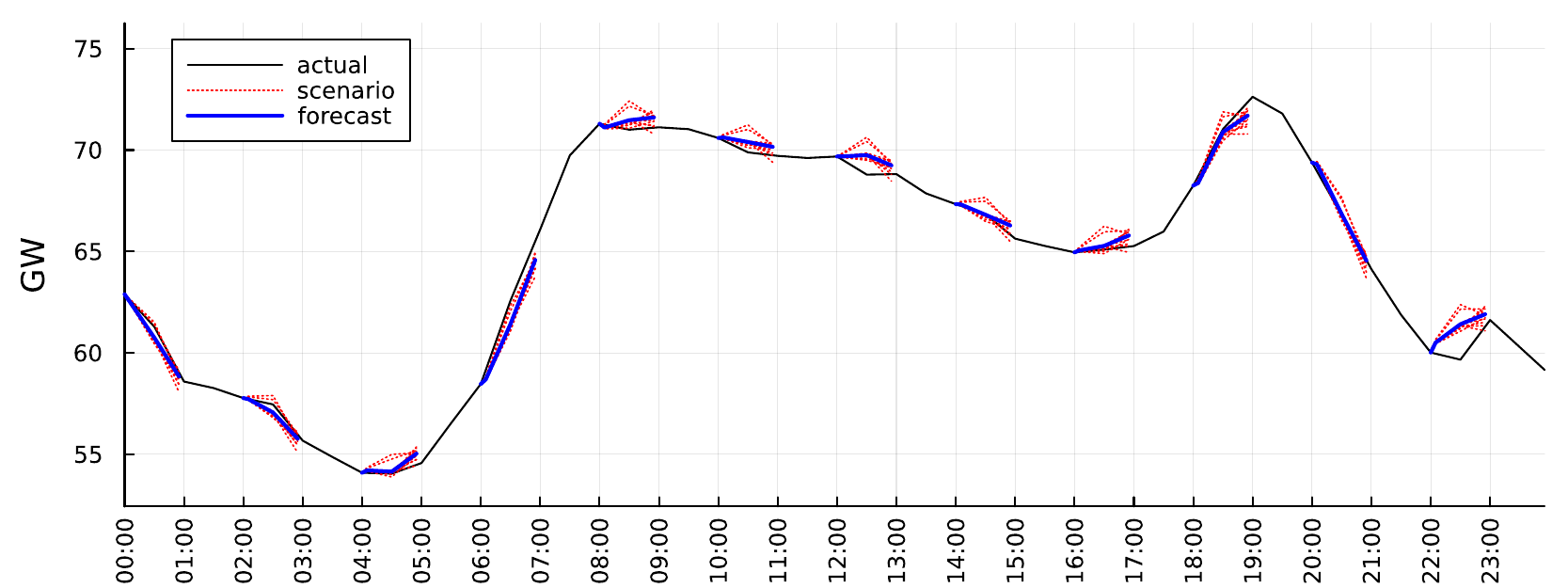}
    \hfill
    \includegraphics[width=0.48\columnwidth]{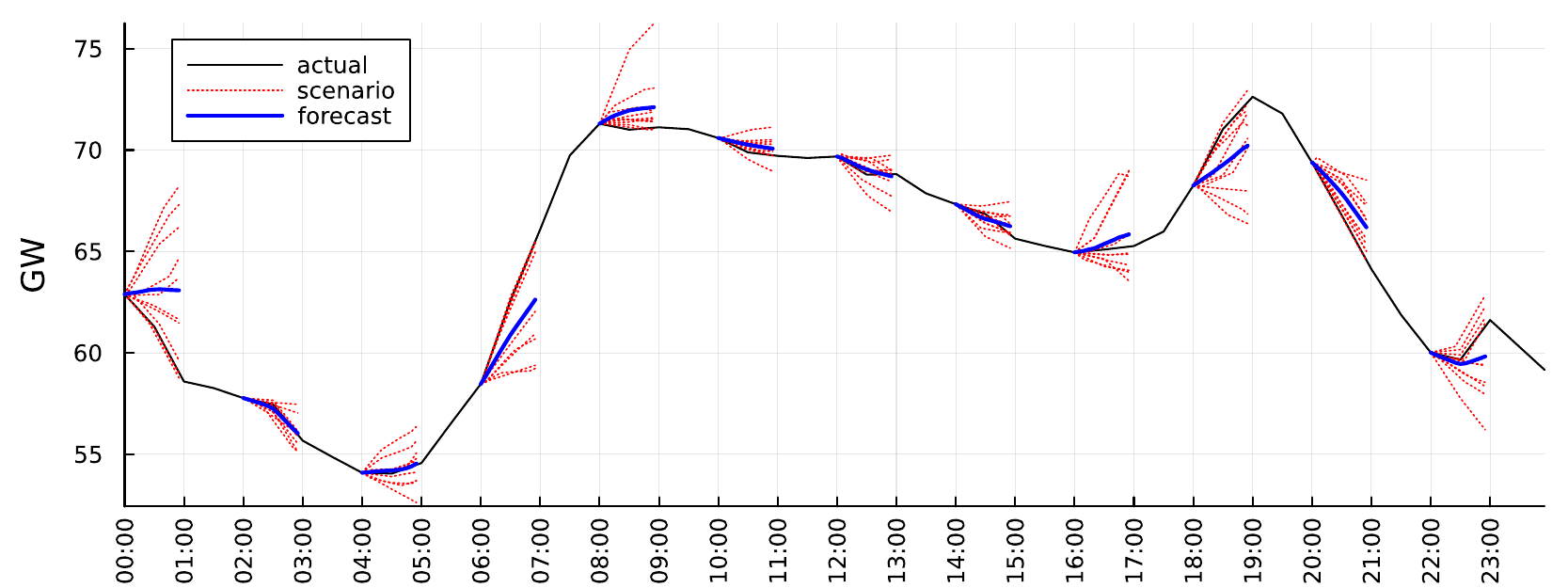}
    \caption{Comparison of probabilistic net load forecasts and scenarios on January 30th, 2018.
        Left: Deep Learning (DL) model.
        Right: K-Nearest Neighbors (KNN) model.
        The DL point prediction is more accurate than KNN, and the DL scenarios have lower variance than KNN.}
    \label{fig:scenarios:2018-01-30:comparison}
\end{figure}

The experiments of Section \ref{sec:experiment_RTE} were also executed with a k-Nearest Neighbor (KNN) method to generate scenarios for LAD and SLAD \citep{cunningham2021k}.
The kNN model selects the $k {=} 10$ closest historical scnearios.
Figure \ref{fig:scenarios:2018-01-30:comparison} illustrates net load probabilistic forecasts produced by the deep learning (DL) models described in Section \ref{sec:experiment_RTE:forecast} on January 30th, 2018.
While the accuracy of the DL model's point forecast is higher than that of the KNN model, the scenarios produced by kNN have higher variance, and systematically cover the actual net load.

Next, Table \ref{tab:total_cost_table:DLvsKNN} reports total operating costs for LAD and SLAD across the year, using both DL and KNN forecasts.
Recall that results obtained with DL forecasts are the same as those reported in Table \ref{tab:total_cost_table:DL}.
First, when using the point forecast produced by KNN, the performance of LAD significantly degrades: on high-load days, LAD+KNN costs $0.8\%$ more compared to LAD+DL.
This highlights the sensitivity of LAD to the quality of the input forecast, i.e., a bad forecast may result in poor decisions.
Indeed, shorter-term forecasts are usually more accurate than longer-term ones.
Therefore, when forecasting errors are factored in, increasing the look-ahead horizon in LAD does not guarantee more cost savings.
In contrast, the performance of SLAD is less sensitive to the input scenarios.
Namely, the year-round average cost of SLAD+KNN is almost identical to that of SLAD+DL, with the former (1.28\%) saving 82.9\% more compared to SCED+RP, and the latter (1.16\%) saving 65.7\% more.
This result suggests that the performance of SLAD is robust to the quality of the point forecast, as long as scenarios capture actual values.

\begin{table*}[!t]
    \centering
    \caption{Comparison of total operating costs of LAD and SLAD averaging over load level in 2018 using different forecasts and scenarios. All costs are in M\$, and all savings are relative to SCED (see Table \ref{tab:total_cost_table:DL}).
    }
    \label{tab:total_cost_table:DLvsKNN}
    \footnotesize
    \begin{tabular}{llrrrrrrrr}
        \toprule
            & \multicolumn{3}{c}{Low load} 
            & \multicolumn{3}{c}{High load}
            & \multicolumn{3}{c}{All}\\
        \cmidrule(lr){2-4}
        \cmidrule(lr){5-7}
        \cmidrule(lr){8-10}
        Problem 
            & Cost
            & Savings
            & Savings(\%)
            & Cost
            & Savings
            & Savings(\%)
            & Cost
            & Savings
            & Savings(\%)\\
        \midrule
        LAD+DL
           & 13.02 & 0.14 & 1.04 \% 
           & 34.47 & 0.36 & 1.04 \% 
           & 23.54 & 0.25 & 1.04 \% 
           \\
        LAD+KNN
           & 13.05 & 0.10 & 0.77 \% 
           & 34.75 & 0.09 & 0.25 \% 
           & 23.69 & 0.09 & 0.40 \% 
           \\
        \midrule
        SLAD+DL
           & \textbf{13.01} & \textbf{0.15} & \textbf{1.11} \% 
           & 34.43 & 0.41 & 1.17 \% 
           & 23.51 & 0.28 & 1.16 \% 
           \\
        SLAD+KNN
           & \textbf{13.01} & \textbf{0.15} & \textbf{1.11} \% 
           & \textbf{34.37} & \textbf{0.47} & \textbf{1.35} \% 
           & \textbf{23.48} & \textbf{0.30} & \textbf{1.28} \% 
           \\
        \bottomrule
    \end{tabular}
\end{table*}

\subsection{Impact of the Look-Ahead Horizon}
\label{sec:sensitivity:horizon}

This section explores the impact of the look-ahead horizon on the performance of SLAD.
Indeed, recall from Table \ref{tab:total_cost_table:DL} that, over the year, SLAD+DL and PD achieve overall savings of $1.16\%$ and $1.55\%$ compared to SCED, respectively.
Also recall that PD cannot be implemented in practice, because it relies on a perfect forecast oracle for the full day.
This motivates considering a longer horizon for SLAD, which may improve its overall economic performance.

To that end, the paper evaluates the performance of SLAD with various look-ahead horizons on January 30th, 2018.
The paper also considers, as a baseline, a LAD formulation with the same horizon as SLAD, but seeded with a perfect forecast; this formulation is referred to as Perfect LAD (PLAD).
It is important to recognize that PLAD formulations are also idealized as they use perfect forecasts and forecasts always incur errors in practice.
Nevertheless, PLAD allows to evaluate the value of a perfect prediction oracle compared to SLAD under the same look-ahead horizons.

\begin{figure}[!t]
    \centering
    \includegraphics[width=0.6\columnwidth]{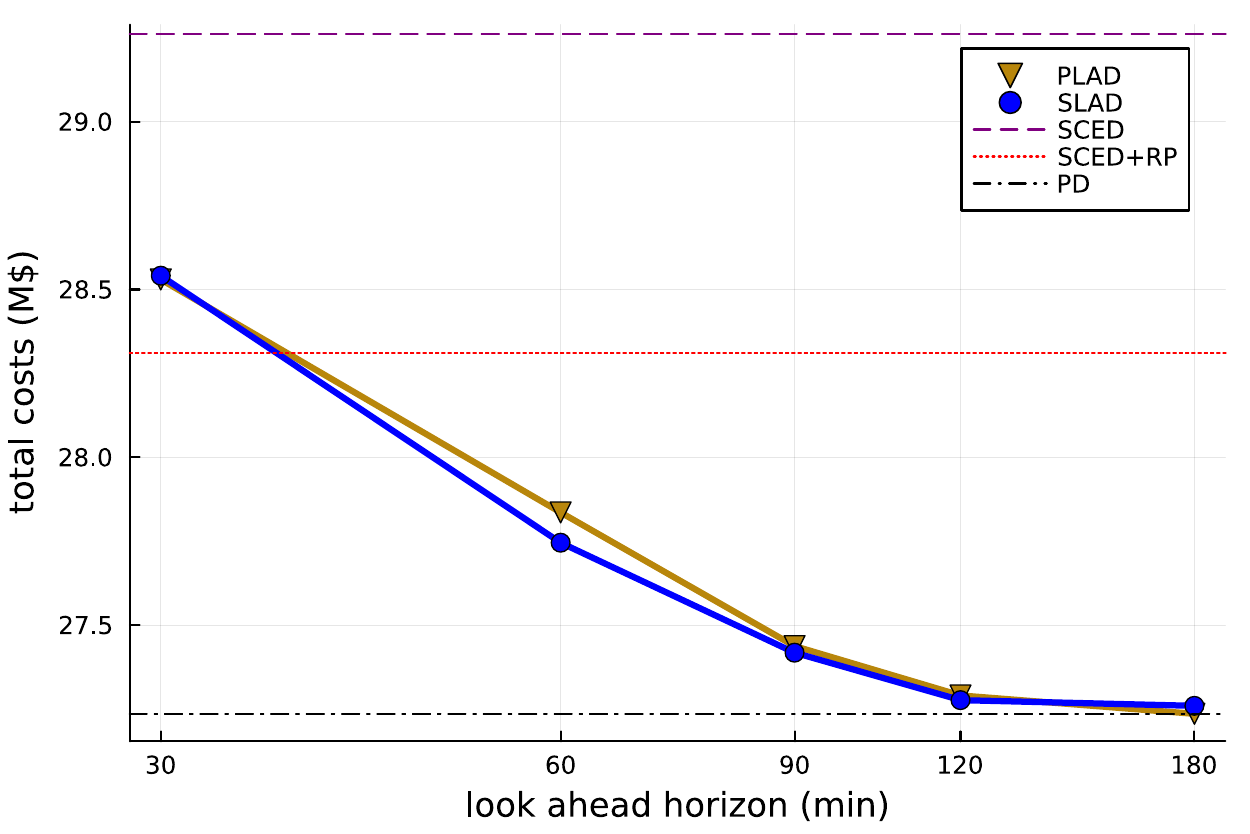}
    \caption{Comparison of SLAD and PLAD with different look-ahead horizons on January 30th, 2018.
    The black dashed line depicts the (unattainable) performance of the idealized Perfect Dispatch (PD).
    The red dashed line depicts the performance of SCED with ramping products (SCED+RP).
    The purple dashed line depicts the performance of SCED.
    }
    \label{fig:2018-01-30:sensitivity:look-ahead-horizon}
\end{figure}

The experiments are conducted with a look-ahead ranging from 30 minutes to 3 hours (180 minutes), using scenarios produced by the KNN model.
Figure \ref{fig:2018-01-30:sensitivity:look-ahead-horizon} illustrates each formulation's total operating costs as a function of the look-ahead horizon.
The figure also depicts the overall costs of PD as an (unreachable) baseline.
First, as expected, the total operational cost of SLAD and PLAD decreases as the look-ahead horizon increases.
Second, there is no significant difference in cost between PLAD with a 3-hour look-ahead and PD.
This result suggests that, for the system at hand, a three-hour horizon is long enough to make nearly optimal dispatch decisions.
Third, the performance of SLAD is almost identical (and sometimes exceeds) to that of PLAD, and almost matches that of PD when a 3-hour horizon is considered.
This demonstrates that the proposed SLAD formulation can produce high-quality solutions without the need for a perfect forecast.

Finally, increasing the look-ahead horizon in SLAD does increase its computational burden, namely, it leads to higher memory requirements and computing time.
Such trade-offs should thus be considered when determining the length of the look-ahead horizon.

\section{Conclusion}
\label{sec:conclusion}

This paper studied the benefits of stochastic look-ahead dispatch in
the context of the MISO's market-clearing pipeline, evaluated on the
France RTE network. The paper showed that multi-period formulations can
yield significant benefits over a single-period formulation. Moreover,
a stochastic formulation (SLAD) brings additional benefits, meets
the market requirements in terms of computing speed, and may even outperform a
multi-period LAD using a perfect forecast for the same horizon.  The
primary benefit of multi-period formulations is their ability to
anticipate ramping events, reducing import costs and transmission line
violations.
The paper also showed that SLAD bridges a significant portion of
the gap between SCED and the potential benefits of a perfect
dispatch over 24 hours. An interesting consequence of these results
concerns the proliferation of ramping products in existing
markets. Multi-period formulations should eliminate, or
at least substantially decrease, the need for such ramping products. 

In terms of future directions, one promising area for improvement is to  further accelerate the SLAD solving process.
Additionally, as renewable generation penetration and storage units continue to increase in power grids, applying SLAD to a system with a larger proportion of renewable generation capacity and storage could provide valuable insights for the future. Another potential direction is to add intra-day unit commitment to the pipeline, which could potentially eliminate some of the import costs highlighted in the paper and lead to a more accurate representation of the proposed stochastic approach. Finally, a systematic comparison of stochastic and robust formulations in terms of computational efficiency, economic performance, and reliability, would be valuable.

\section*{Acknowledgements}

This research was partly supported by NSF award 2112533 and ARPA-E PERFORM award AR0001136.

\bibliographystyle{apalike}
\bibliography{refs.bib}






\end{document}